\documentclass[10pt,final,oneside,onecolumn]{IEEEtran}

\usepackage{graphicx}
\usepackage{amsmath}
\usepackage{amssymb}
\usepackage{wasysym}

\newtheorem{theorem}{Theorem}[section]

\newtheorem{proposition}[theorem]{Proposition}
\newtheorem{lemma}[theorem]{Lemma}

\newtheorem{remark}[theorem]{Remark}

\newcommand{\Hermitian}{{\cal H}(n)}

\newcommand{\Real}{\mathop{\rm Re}}
\newcommand{\Imag}{\mathop{\rm Im}}

\newcommand{\Ker}{\mathop{\rm Ker}}
\newcommand{\tr}{\mathop{\rm tr}}
\newcommand{\Rgamma}{\mathop{\rm Range} \Gamma}
\newcommand{\Jpsi}{{J_\Psi}}
\newcommand{\Lagrange}{L}
\newcommand{\LHG}{{\cal L}^H_\Gamma}
\newcommand{\LH}{{\cal L}^H}
\newcommand{\Hess}{H_\Lambda}

\newcommand{\R}{\mathbb{R}}
\newcommand{\C}{\mathbb{C}}
\newcommand{\N}{\mathbb{N}}
\newcommand{\T}{\mathbb{T}}

\newcommand{\SpecSpace}{{{\cal S}_+^{m\times m}(\mathbb{T})}}
\newcommand{\Spec}{{\cal S}}
\newcommand{\dth}{{\rm d}\vartheta}
\newcommand{\ejth}{{\rm e}^{{\rm j}\vartheta}}
\newcommand{\suchthat}{\mbox{s.\ t.}}
\newcommand{\Unitcircle}{{\mathbb T}}
\newcommand{\qed}{\hfill $\Box$ \vskip 2ex}

\newcommand{\Inversion}{{\cal I}}

\newcommand{\Ball}{B(0, \varepsilon)}

\begin{document}
\title{A globally convergent matricial  algorithm for multivariate spectral estimation}
\author{Federico~Ramponi, Augusto~Ferrante and~Michele~Pavon\thanks{Work partially supported by the MIUR-PRIN Italian grant
``Identification and Control of Industrial Systems"}
\thanks{A. Ferrante and F. Ramponi are  with the
Dipartimento di Ingegneria dell'Informazione, Universit\`a di Padova,
via Gradenigo 6/B, 35131
Padova, Italy {\tt\small  augusto@dei.unipd.it},  {\tt\small  rampo@dei.unipd.it}}\thanks{M. Pavon is with
the Dipartimento di Matematica Pura ed Applicata, Universit\`a di Padova, via
Trieste 63, 35131 Padova, Italy {\tt\small
pavon@math.unipd.it}}}

\markboth{DRAFT}{Shell \MakeLowercase{\textit{et al.}}: Bare Demo of IEEEtran.cls for Journals}

\maketitle

\begin{abstract}
In this paper, we first describe a  {\em matricial} Newton-type algorithm designed
to solve the multivariable spectrum approximation problem. We then prove its 
{\em global} convergence.
Finally, we apply this approximation procedure to {\em multivariate spectral estimation},
and test its effectiveness through simulation.
Simulation shows that, in the case of {\em short observation records}, this method may
provide a valid alternative to standard multivariable identification techniques such as MATLAB's PEM and MATLAB's N4SID.
\end{abstract}


\IEEEpeerreviewmaketitle

\section{Introduction}

ARMA identification methods usually lead to nonconvex optimization problems for which global convergence is not guaranteed, cf. e.g. \cite{LJUNG,SODSTO,STOMOS,Deistler}. Although these algorithms are simple and perform effectively, as observed in \cite[p.103]{STOMOS}, \cite[Section 1]{L}, no theoretically satisfactory approach to ARMA parameter estimation appears to be available. Alternative, convex optimization approaches have been recently proposed by Byrnes, Georgiou, Lindquist and co-workers \cite{BEL2,GL3} in the frame of a broad research effort on analytic interpolation with degree contraint, see
\cite{BLN},
\cite{BEL},
\cite{BGL1},
\cite{BGLM},
\cite{BGuL},
\cite{BGuL2},
\cite{BL2},
\cite{BGLTHREE},
\cite{E},
\cite{GSTATECOV},
\cite{G1},
\cite{G3},
\cite{G4},
\cite{G5},
\cite{G7},
\cite{G8},
\cite{G9},
\cite{G10},
\cite{GLKL},
\cite{GL2}
and references therein.
In particular, \cite{BGLTHREE} describes a new setting for spectral estimation.
This so-called {\em THREE} 
algorithm appears to allow for higher resolution in prescribed frequency bands and to be particularly suitable in case of short observation records. It effectively detects spectral lines and steep variations (see \cite{NEG} for a recent biomedical application). An outline of this method is as follows. A given realization of a stochastic process
(a finite collection of data $y_1 ... y_N$) is fed
to a suitably structured bank of filters, and the steady-state covariance matrix
of the resulting output is estimated by statistical methods. Only zeroth-order covariance lags of the output of the filters need to be estimated, ensuring statistical robustness of the method.
Finding now an input process whose {\em rational} spectrum
is compatible with the estimated covariance
poses naturally  a Nevanlinna-Pick interpolation problem
with bounded degree. The solution of this interpolation problem
is considered as a mean of estimating the spectrum. 
A particular case described in the paper is the {\em maximum differential entropy}
spectrum estimate, which amounts to the so-called central solution
in the Nevanlinna-Pick theory. More generally, the scheme allows for a non constant {\em a priori} estimate $\Psi$ of the spectrum. The Byrnes-Georgiou-Lindquist school has  shown how this and other important problems of control theory may be advantageously 
cast 
in the frame of convex optimization. These problems admit a finite dimensional dual  (multipliers are matrices!) that can be shown to be solvable. The latter result, due to Byrnes and Lindquist \cite{BL2} (see also \cite{FPR2}) is, however, nontrivial since the optimization occurs on an open, unbounded set of Hermitian matrices. The numerical solution of the dual problem is also challenging \cite{BGLTHREE,E,N}, since the gradient of the dual functional tends to infinity at the boundary of the feasible set. Finally, reparametrization of the problem may lead to loss of global concavity, see the discussion in \cite[Section VII]{GLKL}.

This paper adds to this effort in that we consider  estimation of a multivariate spectral 
density in the spirit of THREE \cite{BGLTHREE}, but employing a different metric for the 
optimization part, namely the {\em Hellinger distance} as in \cite{FPRMULTIV}. 
In papers \cite{BGuL, BGuL2}, Byrnes, Gusev and Lindquist chose the 
Kullback-Leibler divergence as a frequency weighted entropy measure,
thus introducing a broad generalization of Burg's maximum entropy method.
More recently, this motivation was supported by
the well-known connection with prediction error methods, 
see e.g. \cite{StoorvogelvanSchuppen,L}.
In the multivariable case, a Kullback-Leibler pseudodistance 
may also be readily defined \cite {G8} inspired by the {\em von Neumann's relative entropy}
 \cite{vonneumann,Ve} of statistical quantum mechanics. The resulting  spectrum approximation 
problem, however, leads to computable solutions of bounded McMillan degree only in the case when 
the prior spectral density is the identity matrix \cite{G8,FPRMULTIV} (maximum entropy solution). 
On the contrary, with a suitable extension of the scalar Hellinger distance introduced in \cite{FPRMULTIV}, 
 the Hellinger approximation generalizes nicely to the multivariable case for any prior estimate $\Psi$ of the spectrum.

The main contributions of this paper, after some background material in Sections II-IV, are found in Sections V-VIII. In Section V, we establish {\em strong convexity} and {\em smoothness}  of the dual functional on a certain domain of Hermitian matrices.
In Section VI, we  analyze in detail a variant of a Newton-type  {\em matricial}  iteration designed to numerically solve the dual of the multivariable spectrum approximation problem. It had originally been sketched in \cite{FPRMULTIV}. The computational burden is dramatically reduced by systematically resorting to solutions of Lyapunov and Riccati equations thanks to various nontrivial results of {\em spectral factorization}. We then show in Section VII that the algorithm is  {\em globally} convergent. 
Finally, in Section VIII, we present guidelines for its application to multivariate spectral estimation and present some simulations comparing to existing methods. Simulation in the multivariable case shows that, at the price of some moderate extra complexity in the model, our method may perform much better than MATLAB's PEM and MATLAB's N4SID in the case of a short observation record.

\section{Constrained spectrum approximation}
Paper \cite{GSTATECOV} introduces and solves the following
moment problem: Given 
a bank of filters described by
an input-to-state stable
transfer function $G(z) = (zI-A)^{-1}B$
and a state covariance matrix $\Sigma$, give necessary and sufficient conditions
for the existence of input spectra $\Phi(\ejth)$ such that
the steady state output has variance $\Sigma$, that is,
\begin{equation}
\label{CONSTRAINT}
\int G\Phi G^\ast = \Sigma.
\end{equation}
Moreover, parametrize the set of all such spectra (here, and in the sequel, integration takes place on the unit circle $\T$ with respect to  normalized Lebesgue measure $d\vartheta/2\pi$) .
Throughout this paper we use the following notations:
$A^\ast = \bar{A}^\top$ for matrices and
$G^\ast \equiv G^\ast(z) = G^\top(z^{-1})$ for spectra
and transfer functions. The scalar product between square matrices is defined
as $\left< A, B \right> = \tr A B^\ast$.
Let $\Spec = \SpecSpace$ be the family of $\C^{m\times m}$-valued
functions defined on the unit circle which are Hermitian,
positive-definite, bounded and coercive.
We have the following {\em existence} result \cite{GSTATECOV}:
There exists $\Phi \in \Spec$ satisfying (\ref{CONSTRAINT})
if and only if there exists $H \in \C^{m \times n}$ such that
\begin{equation}
\label{FEASIBILITY}
\Sigma - A \Sigma A^\ast = BH + H^\ast B^\ast
\end{equation}
Paper \cite{GLKL} deals with the following (scalar) spectrum {\em approximation} problem:
When constraint (\ref{CONSTRAINT}) is feasible, find the spectrum $\Phi$
which minimizes the Kullback-Leibler pseudo distance 
$$d_{KL}(\Psi,\Phi) = \int \Psi \log \frac{\Psi}{\Phi}$$
from an ``a priori'' spectrum $\Psi$, 
subject to the {\em constraint} (\ref{CONSTRAINT}).
It turns out that, if the prior $\Psi$ is rational, the
solution is also rational, and with degree that can be bounded 
in terms of the degrees of $G(z)$ and $\Psi$.
This problem again admits the maximum differential entropy
spectrum (compatible with the constraint)
as a particular case ($\Psi(\ejth) \equiv 1$).
The above minimization poses naturally  a variational problem,
which can be solved using Lagrange theory.
Its dual problem admits a maximum and can be solved exploiting
numerical algorithms.
In \cite{FPRHELL} we restated and solved a similar variational
problem with respect to a different metric, namely the Hellinger distance:
\begin{equation}
\label{SCALARHELLINGER}
d_H(\Psi,\Phi) = \sqrt{\int \left( \sqrt{\Psi} - \sqrt{\Phi} \right)^2 }
\end{equation}
Equation (\ref{SCALARHELLINGER})
defines a {\em bona fide} distance, well-known in mathematical statistics.
The main advantage of this approach to spectral approximation is that
it easily generalizes to the multivariable case,
whereas log-like functionals do not enjoy this property \cite{BLN,G5,G8,FPRMULTIV}. 

\section{Feasibility and the operator $\Gamma$}

\newcommand{\Continuous}{{\mathcal C}(\Unitcircle;{\cal H}(m))}

In this section, we discuss in depth the feasibility of (\ref{CONSTRAINT}). Following \cite{GLKL} and \cite{FPRMULTIV}, let
$\Hermitian = \{M\in \C^{n \times n}: M = M^\ast\}$,
let $\Continuous$ be the
space of ${\cal H}(m)$-valued continuous functions defined on the unit circle,
and let the operator $\Gamma: \Continuous \rightarrow \Hermitian$
be defined as follows:
\begin{equation}
\label{GAMMA}
\Gamma(\Phi) := \int G \Phi G^\ast
\end{equation}
We are interested in the {\em range} of the operator $\Gamma$
which, having to deal with Hermitian matrices, we consider as a 
vector space over the reals.
\begin{proposition} The following facts hold:
\label{RANGEGAMMAPROP} 
\begin{enumerate}
\item Let $\Sigma = \Sigma^\ast > 0$. The following are equivalent:
\begin{itemize}
\item There exists $H \in \C^{m \times n}$ such that identity (\ref{FEASIBILITY}) holds.
\item There exists $\Phi \in \SpecSpace$ such that $\int G \Phi G^\ast = \Sigma$.
\item There exists $\Phi \in \Continuous$, $\Phi > 0$ such that $\Gamma(\Phi) = \Sigma$.
\end{itemize}
\item Let $\Sigma = \Sigma^\ast$ (not necessarily positive definite).
There exists $H \in \C^{m \times n}$ such that identity (\ref{FEASIBILITY}) holds
if and only if $\Sigma \in \Rgamma$.
\item $X \in \Rgamma^\perp$ if and only if $G^\ast(\ejth)X G(\ejth) = 0 \ \forall \vartheta \in [0, 2\pi]$.
\end{enumerate}
\end{proposition}
%

\IEEEproof

As stated above, it was proved in \cite{GSTATECOV} that there
exists $H \in \C^{m \times n}$ such that 
identity (\ref{FEASIBILITY}) holds with Hermitian and positive definite
$\Sigma$ if and only if $\Sigma = \int G \Phi G^\ast$ for some $\Phi \in \SpecSpace$, $\Phi > 0$.
A similar result, albeit with a different algebraic
formulation of the feasibility condition, was proved 
in \cite[Proposition 2.1]{FPRMULTIV}.
The proof of Fact 1 is straightforward once we note that the ``if'' part
of the proof of \cite[Proposition 2.1]{FPRMULTIV} is 
{\em constructive}, and exhibits a {\em continuous}
spectrum. Hence, the fact that there exists a spectrum $\Phi$ such that
$\Sigma = \int G \Phi G^\ast$ is equivalent to there exists a {\em continuous}
spectrum such that the same holds.

As for the second assertion, let $\Sigma \in \Rgamma$. Then
there exists $\Phi \in \Continuous$ such that
\begin{displaymath}
\begin{split}
\Sigma 
&= \int G \Phi G^\ast 
= \int G (\Phi_+ - \Phi_-) G^\ast \\
&= \int G \Phi_+ G^\ast - \int G \Phi_- G^\ast
= \Sigma_+ - \Sigma_-\\
\end{split}
\end{displaymath}
where $\Phi_+$ and $\Phi_-$ are two spectra
such that $\Phi_+ - \Phi_- = \Phi$ (they can be chosen
to be 
bounded away from zero)
and where  
$\Sigma_+$ and $\Sigma_-$
are symmetric positive definite.
Hence $\Sigma$ is a difference of {\em positive} matrices for which (\ref{FEASIBILITY}) holds.
This establishes (\ref{FEASIBILITY}) for $\Sigma$ itself.
Vice versa, suppose that (\ref{FEASIBILITY}) holds for an Hermitian $\Sigma$.
Let $\Sigma_\alpha$ be the unique solution of the following Lyapunov equation:
\begin{displaymath}
\Sigma_\alpha - A\Sigma_\alpha A^\ast = B \left(\alpha B^\ast\right) + \left(\alpha B^\ast\right)^\ast B^\ast = 2\alpha B B^\ast
\end{displaymath}
where $\alpha \in \R$.
Then $\Sigma_\alpha$ depends linearly upon $\alpha$,
i.e.\ $\Sigma_\alpha = \alpha\Sigma_1$, where
$\Sigma_1>0$ since $(A, B)$ is reachable.
Thus, 
there exists an $\alpha$ such that 
$\Sigma_\alpha>0$ and $\Sigma_\alpha > \Sigma$.
Let 
$\Sigma_- = \Sigma_\alpha - \Sigma$.
Then 
$\Sigma_->0$,
and since (\ref{FEASIBILITY}) holds for $\Sigma$ and $\Sigma_\alpha$,
it also holds for 
$\Sigma_-$.
Then assertion 1 implies that 
there exist $\Phi_\alpha>0$ and $\Phi_2>0$ in $\Continuous$ such that
$\Sigma_\alpha = \int G \Phi_\alpha G^\ast$ and
$\Sigma_- = \int G \Phi_2 G^\ast$,
hence $\Sigma = \int G (\Phi_\alpha - \Phi_2) G^\ast$
and assertion 2 follows. 

The third assertion is a simple geometrical fact: If $X \in \Rgamma^\perp$, then
for any $\Phi \in \Continuous$
\begin{displaymath}
0 = \left<X, \int G\Phi G^\ast \right>
 = \tr X \int G\Phi G^\ast
 = \tr \int (G^\ast X G) \Phi
\end{displaymath}
and the conclusion follows.
\qed

\newcommand{\Span}{\mathop{\rm span}}
\newcommand{\spanSpec}{{\Spec_2}}
\newcommand{\spanGamma}{{\Gamma_2}}
\begin{remark}
The underlying statement in Proposition \ref{RANGEGAMMAPROP},
Facts 1 and 2, is that if we defined $\Gamma$ over 
the vector space of finite linear combinations of functions
in $\SpecSpace$, its range would remain the same.
\end{remark}

\begin{remark}
\label{RANGEGAMMAREMARK} 
Proposition \ref{RANGEGAMMAPROP} shows that $\Rgamma$ is the set of all the 
Hermitian matrices $\Sigma$ for which there exists $H$ such that 
(\ref{FEASIBILITY}) holds. 
This fact will be useful in numerical computations. 
Indeed, $\Rgamma$ is obviously finite-dimensional,
and if $\{H_1, ..., H_{N}\}$ is a base
of $\C^{m\times n}$, then the corresponding solutions
$\{\Sigma_1, ..., \Sigma_{N}\}$ of (\ref{FEASIBILITY}),
{\em considered as a discrete-time Lyapunov equation in the unknown
$\Sigma$}, generate $\Rgamma$.
Note that $\{\Sigma_1, ..., \Sigma_{N}\}$
are not necessarily linearly independent.
\end{remark}

\section{Multivariable spectrum approximation in the Hellinger distance}
Let the function $d_H: \SpecSpace \times \SpecSpace \rightarrow \R^+$
be defined as follows:
\begin{equation}
\label{HELLINGER}
\begin{split}
d_H(\Psi,\Phi)^2 &:= \inf_{W_\Psi, W_\Phi} \tr \int \left(W_\Psi-W_\Phi\right)\left(W_\Psi-W_\Phi\right)^\ast \\
&\suchthat \quad W_\Psi W_\Psi^\ast = \Psi, \quad W_\Phi W_\Phi^\ast = \Phi, \\
&\quad\quad\ W_\Psi\ {\rm and}\ W_\Phi\ {\rm of\ dimension}\ m\times m
\end{split}
\end{equation}
that is, $d_H(\Psi,\Phi)$ is the $L^2$ distance between the sets of
{\em square spectral factors} of the two spectra.
It was shown in \cite{FPRMULTIV} that the infimum in (\ref{HELLINGER}) 
is actually a minimum and that $d_H$ is a {\em bona fide}
distance between spectral densities, and reduces to the ordinary Hellinger distance
in the scalar case. It was also shown there that the minimum in (\ref{HELLINGER}) 
is the same if we fix a square spectral factor $W_\Psi$ of $\Psi$,
and then minimize over the spectral factors of $\Phi$:
\begin{equation}
\label{HELLINGER2}
\begin{split}
d_H(\Psi,\Phi)^2 &\equiv \min_{W} \tr \int \left(W_\Psi-W\right)\left(W_\Psi-W\right)^\ast \\
&\suchthat \quad W W^\ast = \Phi
\end{split}
\end{equation}
The multivariable spectrum approximation problem addressed in \cite{FPRMULTIV} is the following.
Let $G(z) = (zI-A)^{-1}B$, where $A$ is stable, $B$ has full rank and $(A,B)$ is reachable.
Given $\Psi \in \Spec$, find
\begin{equation}
\label{PRIMALPROBLEM}
\arg\min_{\Phi \in \Spec} d_H(\Psi,\Phi) \ \suchthat \ \int G\Phi G^\ast = \Sigma
\end{equation}

Since $\Sigma>0$,
applying the following change of base to $G(z)$:
\begin{displaymath}
\begin{split}
\bar{A}    &= \Sigma^{-1/2} A \Sigma^{1/2}, \\
\bar{B}    &= \Sigma^{-1/2} B, \\
\bar{G}(z) &= (zI - \bar{A})^{-1}\bar{B} = \Sigma^{-1/2} G(z),
\end{split}
\end{displaymath}
it is easy to see that there
is no loss of generality in taking $\Sigma = I$.
Thus, once a factor $W_\Psi$ of $\Psi$ is fixed,
(\ref{PRIMALPROBLEM}) reduces to:
\begin{equation}
\label{PRIMALPROBLEM2}
\begin{split}
&\arg\min_{W} \tr \int \left(W_\Psi-W\right)\left(W_\Psi-W\right)^\ast \\
&\suchthat \quad \int G W W^\ast G^\ast = I
\end{split}
\end{equation}
which is an $L^2$ constrained minimization.

\begin{remark}
In \cite[Theorem 6.1]{FPRMULTIV}, it was shown that the minimizer in (\ref{HELLINGER2}) is explicitly given by
\begin{equation}
\label{PHIHAT}
\hat{W}_\Phi = \Phi^{1/2}[\Phi^{1/2}\Psi\Phi^{1/2}]^{-1/2}\Phi^{1/2}W_\Psi.
\end{equation}
Nevertheless, to solve the approximation problem 
(\ref{PRIMALPROBLEM2}), 
we do not need to employ (\ref{PHIHAT}) (see below). 
\end{remark}

\subsection{Variational analysis}

Now let us assume that the problem is feasible, i.e., condition (\ref{FEASIBILITY}) holds.
To solve (\ref{PRIMALPROBLEM2}), form the Lagrangian:
\begin{equation}
\label{LAGRANGIAN}
\begin{split}
\Lagrange(W, \Lambda) &= \tr \int \left(W_\Psi-W\right)\left(W_\Psi-W\right)^\ast\\ 
&+ \left<\Lambda, \int G W W^\ast G^\ast - I\right>
\end{split}
\end{equation}
where
$\Lambda \in \Hermitian$.
Since $\int G W W^\ast G^\ast \in \Rgamma$ by construction, and $I \in \Rgamma$
{\em by the feasibility assumption}, it is natural,
though not strictly necessary, to restrict {\em a priori}
the Lagrange parameter $\Lambda$ to $\Rgamma$ (a $\Lambda \in \Rgamma^\perp$
would not play any role in the above Lagrangian).

We proceed with {\em unconstrained} minimization of (\ref{LAGRANGIAN}).
The functional (\ref{LAGRANGIAN}) is convex and differentiable in $W$.
Thus, to find the unique minimizing solution
we impose that the first variation of (\ref{LAGRANGIAN})
is zero in each direction $\delta W$.
We easily find the following condition for $W$
(see \cite{FPRMULTIV} for the details):
\begin{displaymath}
W - W_\Psi + G^\ast \Lambda G W = 0
\end{displaymath}
To carry on with the computations, and to ensure that the resulting 
optimum spectrum is integrable over the unit circle, we require
{\em a posteriori} that $\Lambda$ belongs to the following set:
\begin{equation}
\label{LH}
\LH = \left\{ \Lambda \in \Hermitian, \ \ I + G^\ast \Lambda G > 0 \ \ \forall\ \ejth \in \Unitcircle \right\}
\end{equation}
that is, $\Lambda \in \LHG$, where
\begin{equation}
\label{LHG}
\LHG := \Rgamma \cap \LH.
\end{equation}
If this is the case, the optimal spectral factor and the
corresponding optimal spectral density are easily found to be:
\begin{equation}
\label{PRIMALMINIMUM}
\begin{split}
\hat{W}    &= (I + G^\ast \Lambda G)^{-1} W_\Psi, \\
\hat{\Phi} &= (I + G^\ast \Lambda G)^{-1} \Psi (I + G^\ast \Lambda G)^{-1}
\end{split}
\end{equation}
\begin{remark}Observe that, when $\Psi$ is rational, (\ref{PRIMALMINIMUM}) yields a  rational spectrum with McMillan degree that can be bounded. The same applies to scalar spectrum approximation problem in a Kullback-Leibler type distance where the degree of the optimal approximant is actually lower \cite{GLKL}. In the multivariable case, however, the Kullback-Leibler solution is computable and of bounded McMillan degree only when $\Psi=I$ (maximum entropy solution), see  \cite[Theorem 1]{G5}, \cite[Section 4]{G8}. 
\end{remark}

Consider now the issue of {\em existence} of a matrix $\Lambda \in \LHG$
such that 
\begin{equation}
\label{MINIMUMSPECSATISFIESCONSTRAINT}
\int G (I + G^\ast \Lambda G)^{-1} \Psi (I + G^\ast \Lambda G)^{-1} G^\ast = I
\end{equation}
that is, such that the corresponding optimal spectrum satisfies the
constraint (\ref{CONSTRAINT}).
A key result of \cite[Theorem 7.7]{FPRMULTIV},
inspired by a fundamental result of Byrnes and Lindquist \cite{BL2},
states that such a $\Lambda$ indeed
exists and is unique, therefore establishing that for such $\Lambda$
(\ref{PRIMALMINIMUM}) is the solution
to problem (\ref{PRIMALPROBLEM2}).

\begin{remark}
Identity (\ref{MINIMUMSPECSATISFIESCONSTRAINT}) is attained
at $\Lambda = 0$ if and only if $\Psi$ itself satisfies the contraint
(\ref{CONSTRAINT}). The ``only if'' part is trivial.
As for the ``if'' part, let $\int G\Psi G^\ast = I$. Then substituting
$\Lambda = 0$ into (\ref{PRIMALMINIMUM}) we obtain the spectrum
$\hat{\Phi} \equiv \Psi$,
which has the least possible distance $d_H = 0$ from $\Psi$ and
hence is automatically optimal, and
trivially satisfies (\ref{MINIMUMSPECSATISFIESCONSTRAINT}).
The assertion follows from the uniqueness of $\Lambda$.
\end{remark}

In order to find the optimal $\Lambda$, we form the dual functional
(see \cite{FPRMULTIV}):
\begin{equation}
\label{DUALFUNCTIONAL}
\Lagrange_d(\Lambda) = \Lagrange(\hat{W}, \Lambda)  
= \tr \int \left( \Psi - (I + G^\ast \Lambda G)^{-1} \Psi \right) - \tr \Lambda 
\end{equation}
Notice that $\Lagrange_d(\Lambda)$ is finite on $\LHG$.
Recall that the dual of a Lagrangian functional is always {\em concave}
and that a finite convex (or concave) function
defined on a finite dimensional space is continuous on the interior of its
domain (see \cite{KOSMOL} or \cite{ROCKAFELLAR}).
Instead of maximizing (\ref{DUALFUNCTIONAL}), we consider the equivalent
problem of minimizing the following functional:
\begin{equation}
\label{JPSI}
\Jpsi(\Lambda) = - \Lagrange_d(\Lambda) + \tr\int\Psi
= \tr \int (I + G^\ast \Lambda G)^{-1} \Psi + \tr \Lambda 
\end{equation}
The minimization of the
convex and continuous functional 
$\Jpsi$ over $\LHG$ is the main subject of this paper.
The following sections are dedicated to prove strict convexity and smoothness of $\Jpsi$,
to describe a Newton-type algorithm for its numerical minimization,
and to prove the global convergence of that algorithm. Some numerical simulations
follow.

\section{Properties of the functional $J_\Psi$}
\newcommand{\Ql}{Q_\Lambda}
\newcommand{\Qinv}{Q_\Lambda^{-1}}
\newcommand{\Qli}{Q_{\Lambda_i}}
\newcommand{\Qiinv}{Q_{\Lambda_i}^{-1}}

\newcommand{\SmallL}{\delta\Lambda}
\newcommand{\SmallMatrix}{M}
\newcommand{\SmallSequence}{\SmallMatrix_n}

\newcommand{\gradient}[2]{{\nabla J_{\Psi,{#1}}(#2)}}

{In this section}, we establish various properties of the functional $\Jpsi$ on $\LHG$.
We begin by recalling a few basic definitions and facts 
from multivariate analysis.
A function $f:S\subset \R^N \rightarrow \R^M$ is ({\em Fr\'echet}) differentiable on
the open set
$S$ if for all $x\in S$ there exists a linear map
$L_x: \R^N \rightarrow \R^M$ such that
\begin{displaymath}
\lim_{h\rightarrow 0}
\frac{||f(x+h) - f(x) - L_x(h)||}{||h||} = 0.
\end{displaymath}
A function
$f$ is said to be $C^0(S)$ if it is continuous on $S$.
Also,
$f$ is said to be $C^1(S)$ if it is differentiable at each $x\in S$
and if the operator $Df$ defined by
\begin{displaymath}
Df(x) := L_x
\end{displaymath}
is $C^0(S)$.
Now 
the derivative
$Df:S\rightarrow L(\R^N, \R^M) \simeq \R^{MN}$
is itself a function between finite-dimensional spaces.
If $Df$ is $C^1(S)$, 
then
$f$ is said to be $C^2(S)$. 
Proceeding in this way, the $C^k$-differentiability of $f$
can be defined. 
Finally, $f$ is said to be of class $C^\infty(S)$ if it is $C^k(S)$ for all $k\in \N$.
A standard result in analysis 
states that $f\in C^1(S)$ if and only if
the partial derivatives $\frac{\partial f_m}{\partial x_n}$
(where $f_m$ is the $m$-th component of $f$) 
exist and are continuous on $S$ 
(see for instance \cite[Theorem 9.21]{RudinPMA}).
It follows that $f\in C^k(S)$ if and only if $f$ has in $S$ 
continuous partial derivatives of any order up to $k$,
that is:
\begin{displaymath}
\frac{\partial^h f_m}
{\partial x_{n_1} \cdots \partial x_{n_h}}
\in C^0(S) 
\end{displaymath}
for all $m, n_i, h$ s.t.\ 
$1 \leq m \leq M$,
$1 \leq n_i \leq N$,
and $0 \leq h \leq k$.

In our setting, where $\Jpsi: \LHG\subset\Rgamma \rightarrow \R$,
the role of the above partial derivatives is played
by the directional (or {\em G\^ateaux}) derivatives,
\begin{displaymath} 
\delta^h \Jpsi(\Lambda; \bar{\delta\Lambda}_{n_1}, ..., \bar{\delta\Lambda}_{n_h})
\end{displaymath}
where $1 \leq n_i \leq d$ and
$\{\bar{\delta\Lambda}_1, ..., \bar{\delta\Lambda}_d\}$ is a fixed
orthonormal base of $\Rgamma$.
{\em A fortiori}, if we show that $\Jpsi$ has on $\LHG$ {\em continuous} directional 
derivatives of any order up to $k$, taken in {\em whatever} directions
$\{\delta\Lambda_1, ..., \delta\Lambda_k\} \subset \Rgamma$,
then we can say that $\Jpsi \in C^k(\LHG)$.

\begin{lemma}
\label{LEMMA0}
Let $H\in\Hermitian$ and $m$ be its minimum eigenvalue.
The map $H \mapsto m$ is continuous. 
\end{lemma}

\IEEEproof

The map from a matrix $H$ to the vector of coefficients
of its characteristic polynomial $a(s) = \det(sI - H) = a_0 + ... + a_{n-1}s^{n-1} + s^n$ is continuous.
Indeed, each of the coefficients of $a(s)$ is obtained by means of sums and products of  elements of $H$.
Moreover,
it is a well-known fact (see for example \cite{MARDEN}) 
that the mapping from the coefficients
of a monic polynomial to its roots is continuous, in the following
sense: Given
$a(s)=s^n+\sum_{i=0}^{n-1}a_is^i$, let $\lambda_i$ be the zeros of $a(s)$ and $\nu_i$ the respective multiplicities. 
For all $\varepsilon > 0$, there exists $\delta > 0$ such that
if $b(s) = s^n+\sum_{i=0}^{n-1}b_is^i$ and
$|b_i - a_i|<\delta$ for all $i=0,1,\dots,n-1$, then $b(s)$ has $\nu_i$ zeros in the ball
centered in $\lambda_i$ with radius $\varepsilon$.
In conclusion, if $H$ is Hermitian, the mapping from $H$ to its minimum
(real) eigenvalue is continuous. 
\qed

\begin{lemma}
\label{Lemma} 
Define $\Ql(z) = I + G^\ast(z) \Lambda G(z)$.
Consider a sequence $\Lambda_n\in\LHG$ converging to $\Lambda\in\LHG$. 
Then $Q^{-1}_{\Lambda_n}$ are  well defined and continuous on $\T$
and converge {\em uniformly} to $\Ql^{-1}$ on $\T$.
\end{lemma}
\noindent
\IEEEproof

Observe that, for  $\Lambda \in \LHG$, $\Ql$ is a positive definite, 
continuous matrix function on $\T$. 
By Lemma \ref{LEMMA0}, there exists a continuous function $m_\Lambda(\ejth)>0$ such 
that $\Ql(\ejth)\ge m_\Lambda(\ejth) I$ for each 
$\vartheta$. 
Hence, $\Ql(\ejth)\ge m_\Lambda I, \forall \vartheta$, where 
$m_\Lambda:=\min_\vartheta m_\Lambda(\ejth)>0$. Let $\SmallL \in \Ball$, 
the {\em closed} ball of radius $\varepsilon$ centered in $0$.
Now,  
\begin{displaymath}
||G^\ast(\ejth)\SmallL G(\ejth)|| \leq ||\SmallL|| M_G \leq \varepsilon M_G
\end{displaymath}
where
\begin{displaymath}
M_G = \max_\vartheta ||G^\ast(\ejth)|| \ ||G(\ejth)||.
\end{displaymath}
Thus, if we choose $\varepsilon < m_\Lambda / M_G$,
then $||G^\ast(\ejth)\SmallL G(\ejth)|| < m_\Lambda$. Hence,
$I + G^\ast(\Lambda+\SmallL)G$ describes, as $(\SmallL, \vartheta)$
varies in $\Ball \times [-\pi,\pi]$,
a compact set that does not contain any singular matrix.
Now recall that the matrix inversion operator is continuous at any nonsingular
matrix. Hence, $Q_{\Lambda+\SmallL}^{-1} (\ejth)$
admits a uniform bound $M(\Lambda,\epsilon)$ on $\Ball \times [-\pi,\pi]$.
Since $\Lambda_n\rightarrow \Lambda$, for $n$ sufficiently large,
$(\Lambda_n-\Lambda)\in\Ball$. Then
\begin{displaymath}
\begin{split}
\sup_\vartheta\|Q^{-1}_{\Lambda_n}-\Ql^{-1}\|
  &=       \sup_\vartheta \|Q^{-1}_{\Lambda_n} \left[ G^\ast(\Lambda-\Lambda_n)G \right] \Ql^{-1}\| \\
  &\le M^2 \sup_\vartheta \|G^\ast(\Lambda-\Lambda_n)G \| \\
  &\le M^2 \varepsilon M_G.
\end{split}
\end{displaymath}
This implies that $Q^{-1}_{\Lambda_n}\rightarrow\Ql^{-1}$ uniformly on $\T$.
\qed

\begin{theorem} \label{JPSIREGULARITYPROP} Consider $\Jpsi: \LHG\subset\Rgamma \rightarrow \R$. Then
\begin{enumerate}
\item $\Jpsi \in C^\infty(\LHG)$.
\item $\Jpsi$ is strictly convex on $\LHG$.
\end{enumerate}
\end{theorem}

\noindent
\IEEEproof
 
Let $\Inversion: A \mapsto A^{-1}$ be the matrix inversion operator.
Making use of 
\begin{equation}\label{INVERSIONDERIVATIVE}
\delta\Inversion (A; \delta A) = - A^{-1} \delta A A^{-1}, 
\end{equation}
the first variation of $\Jpsi(\Lambda)$ in an arbitrary direction
$\delta\Lambda_1$ is found to be:
\begin{equation}
\label{JPSIFIRST}
\begin{split}
\delta\Jpsi(\Lambda; \delta\Lambda_1) 
&= - \tr \int \Qinv G^\ast \delta\Lambda_1 G \Qinv \Psi + \tr \delta\Lambda_1 \\
&= \left< I - \int G \Qinv \Psi \Qinv G^\ast, \quad \delta\Lambda_1 \right> 
\end{split}
\end{equation}
The linear functional 
$\gradient{\Lambda}{\cdot} := \delta\Jpsi(\Lambda; \cdot)$ 
defined by (\ref{JPSIFIRST}) is the {\em gradient} of $\Jpsi$ at $\Lambda$.
To prove that $\Jpsi \in C^1(\LHG)$
we must show that, for any fixed $\delta\Lambda_1$, 
$\delta\Jpsi(\Lambda; \delta\Lambda_1)$
is continuous in the variable $\Lambda$
(it follows that 
$\gradient{\Lambda}{\cdot}$
is also continuous in $\Lambda$).
Consider a sequence $\SmallSequence\in\Rgamma$ converging to $0$. 
By Lemma \ref{Lemma},
$Q_{\Lambda+\SmallSequence}^{-1}$ converge uniformly to
$\Ql^{-1}$. 
Recall that $\Psi$ 
is bounded. Applying elementwise the bounded convergence theorem, we get
\begin{displaymath}
\lim_{n\rightarrow + \infty}
\int Q_{\Lambda+\SmallSequence}^{-1} G^\ast \delta\Lambda_1 G Q_{\Lambda+\SmallSequence}^{-1} \Psi
= \int \Qinv G^\ast \delta\Lambda_1 G \Qinv \Psi
\end{displaymath}
Hence, for all $\delta\Lambda_1 \in \Rgamma$,
$\delta\Jpsi(\Lambda; \delta\Lambda_1)$
is continuous, i.e.\ $\Jpsi \in C^1(\LHG)$.
The second variation of $\Jpsi$, say in direction $\delta\Lambda_2$, is easily
obtained applying (\ref{INVERSIONDERIVATIVE}) and the chain rule
to (\ref{JPSIFIRST}):
\begin{equation}
\label{JPSISECOND}
\begin{split}
&\delta^2\Jpsi(\Lambda; \delta\Lambda_1, \delta\Lambda_2) \\
&= \tr \int W_\Psi^\ast \Qinv G^\ast \delta\Lambda_2 G \Qinv G^\ast \delta\Lambda_1 G \Qinv W_\Psi  \\
&+ \tr \int W_\Psi^\ast \Qinv G^\ast \delta\Lambda_1 G \Qinv G^\ast \delta\Lambda_2 G \Qinv W_\Psi  
\end{split}
\end{equation}
The bilinear form $H_\Lambda(\cdot, \cdot) := \delta^2\Jpsi(\Lambda; \cdot, \cdot)$
is the {\em Hessian} of $\Jpsi$ at $\Lambda$.
Again, continuity of 
$\delta^2\Jpsi(\Lambda; \delta\Lambda_1, \delta\Lambda_2)$
can be established by the previous argument in view of Lemma \ref{Lemma}.
Similarly, it can be shown that $\Jpsi$ has continuous directional
derivatives of any order. Thus $\Jpsi \in C^k(\LHG)$ for any $k$, and the first assertion follows.
Finally, we show that $\Jpsi$ is {\em strictly} convex on $\LHG$.
A standard result in the theory of convex functions states
that if a function $f: S\subset \R^N \rightarrow \R$ is
$C^2(S)$ (where $S$ is open), then $f$ is
{\em strictly convex} on $S$ if and only if its Hessian $H_x$
is positive definite at each $x\in S$. Consider $H_\Lambda(\delta\Lambda, \delta\Lambda)=\delta^2\Jpsi(\Lambda; \delta\Lambda, \delta\Lambda)$ for $\delta\Lambda \in \Rgamma \setminus \{0\}$. Since the integrand in (\ref{JPSISECOND}) is   positive semidefinite, it follows that $H_\Lambda(\delta\Lambda, \delta\Lambda)\ge 0$.
In view of Point 3 in Proposition \ref{RANGEGAMMAPROP}, the integrand is not identically zero and $H_\Lambda(\delta\Lambda, \delta\Lambda)> 0$.
It follows that $\Jpsi$ is strictly convex.
\qed

Loosely speaking,
\cite{FPRMULTIV} establishes existence of the minimum by showing
that $\Jpsi$ is bounded from below and that
$\Jpsi(\Lambda) \rightarrow +\infty$ whenever
$||\Lambda|| \rightarrow +\infty$ or $\Lambda$ approaches
$\partial\LHG$, the {\em boundary} of $\LHG$. Since
the minimum point cannot reside at infinity or at the boundary,
we can then restrict the minimization problem to a sublevel set
of $\Jpsi$. It follows from the continuity of $\Jpsi$ that
such a set is compact, and the existence result follows
from Weierstrass' theorem.
Since $\Jpsi$ is strictly convex, the minimum point is {\em unique}. 

\section{A matricial Newton algorithm}
\subsection{Description of the iterative method}

The Newton algorithm is an iterative procedure for the
search of roots of a function or the {\em minimization}
of a functional.
With respect to the latter objective, it can be formulated as follows.
Let $f:S\rightarrow \R$ be a functional defined over $S\subset \R^n$.
In order to find an estimate $\hat{x}$ of a minimum point $x^\ast$ of $f$,
\begin{enumerate}
\item Make an initial guess $x_0$, possibly near the minimum point.
\item At each iteration, compute the Newton step
\begin{equation}
\Delta x_i = -H_{x_i}^{-1}\nabla f_{x_i}
\end{equation}
where $H_{x_i}$ is the Hessian of $f$ at $x_i$ and $\nabla f_{x_i}$ is
the gradient of $f$ at $x_i$
(understood as a column vector).
\item Set $t_i^0 = 1$, and let $t_i^{k+1} = t_i^k/2$
until both of the following conditions hold:
\begin{equation}
\label{BACKTRACKINGDOMAIN}
x_{i} + t_i^k\Delta x_i \in S
\end{equation}
\begin{equation}
\label{BACKTRACKING}
f(x_i + t_i^k\Delta x_i) < f(x_i) + \alpha t_i^k \nabla f_{x_i}^\top \Delta x_i
\end{equation}
where $\alpha$ is a real constant, $0<\alpha<1/2$.
\item Set $x_{i+1} = x_{i} + t_i^k\Delta x_i$.
\item Repeat steps 2, 3 and 4 until $|\nabla f_{x_i}| < \varepsilon$,
where $\varepsilon$ is a (small) tolerance threshold, then set
$\hat{x} = x_i$.
\end{enumerate}
In its ``pure'' form, the iteration of the Newton algorithm
only consists in step 2, which is indeed its essential part.
Step 3 is the so-called {\em backtracking} procedure.
For small $t$, if $f$ is sufficiently regular, we have
$f(x_i + t\Delta x_i) \simeq f(x_i) + t \nabla f_{x_i}^\top \Delta x_i$.
Since $\nabla f_{x_i}^\top \Delta x_i = - \nabla f_{x_i}^\top H_{x_i}^{-1}\nabla f_{x_i} < 0$,
condition (\ref{BACKTRACKING}) {\em must} hold for small $t$,
hence step 3 must terminate at some iteration. 
Since $\nabla f_{x_i}^\top \Delta x_i < 0$,
(\ref{BACKTRACKING}) implies $f(x_i + t_i^k\Delta x_i) < f(x_i)$.
That is, $\{f(x_i)\}$ is a strictly decreasing sequence.

In essence, the ``pure'' Newton algorithm works {\em very} well when
the starting point happens to be near the minimum and the
function $f$ is there effectively approximated by a quadratic form,
but it can suffer from numerical problems when this is not the case.
The backtracking line search is a remedy to this drawback;
moreover it can be shown that, under certain regularity assuptions on $f$,
{\em which hold in our case} (see Section \ref{SECTIONCONVERGENCE}),
after a finite number of iterations step 3 always selects the
multiplier $t = 1$, that is, the full step.
During the latter stage, the convergence to the minimizing solution
is {\em quadratic}, meaning that there exists a constant $C$
such that $||x_{i+1} - x^\ast|| \leq C||x_i - x^\ast||^2$.
This rate of convergence makes the Newton algorithm often
preferable over other minimization methods
(see \cite{BV}).
We must minimize the functional $\Jpsi(\Lambda)$ over the set
$\LHG$.
As initial condition, we can safely choose $0$. Hence, set
\begin{equation}
\label{NEWTONSTEP0}
 \Lambda_0 = 0.
\end{equation}
It turns out that, although the problem
is finite-dimensional, the {\em inversion of the Hessian} is more
demanding than inverting  a matrix.
In order to compute the Newton step $\Delta\Lambda_i$,
we must solve at $\Lambda_i$ the following linear equation:
\begin{equation}
\label{NEWTONSTEP}
H_{\Lambda_i}(\Delta\Lambda_i, \cdot) = - \gradient{\Lambda_i}{\cdot}
\end{equation}
where, once fixed $\Lambda_i$,
$\gradient{\Lambda_i}{\cdot}$ 
and $H_{\Lambda_i}(\cdot, \cdot)$ must
be understood as a linear and a bilinear form,
defined
by (\ref{JPSIFIRST}) and (\ref{JPSISECOND}) respectively.
Comparing with the above definitions, (\ref{NEWTONSTEP}) reduces to:
\begin{equation}
\label{NEWTONSTEP2}
\begin{split}
&   \int G \Qiinv \left[ (G^\ast (\Delta\Lambda_i) G \Qiinv \Psi) + (G^\ast (\Delta\Lambda_i) G \Qiinv \Psi)^\ast \right] \Qiinv G^\ast  \\  
& = \int G \Qiinv \Psi \Qiinv G^\ast - I 
\end{split}
\end{equation}
In principle, equation (\ref{NEWTONSTEP2}) is not difficult to solve.
We suggest the following procedure:
\begin{itemize}
\item At the beginning of the procedure,
take a base $\{H_1, ..., H_k, ..., H_N\}$ of  
$\C^{m\times n}$.\ \footnote{
Actually, it suffices to take the $\{H_k\}$ to be
a base of $\C^{m\times n} \ominus \Ker R$, where the map 
$R$ is defined by $$R: H \mapsto BH + H^\ast B^\ast.$$
}
Then compute the solutions 
$\{\Sigma_1, ..., \Sigma_k, ..., \Sigma_N\}$ of the
following discrete-time Lyapunov 
equations:
\begin{displaymath}
\Sigma_k - A \Sigma_k A^\ast = BH_k + H_k^\ast B^\ast
\end{displaymath}
As shown before, these solutions generate $\Rgamma$.

\item To compute $\Delta\Lambda_i$ at each step,
\begin{enumerate} 
\item Compute the integral
\begin{equation}
\label{INTEGRAL1}
Y = \int G \Qiinv \Psi \Qiinv G^\ast - I
\end{equation}
\item For each $\Sigma_k$ in the precomputed generators, compute the following integral:
\begin{equation}
\label{INTEGRAL2}
\begin{split}
Y_k = \int & G \Qiinv \left[ (G^\ast \Sigma_k G \Qiinv \Psi) \right. \\
           &\left.+ (G^\ast \Sigma_k G \Qiinv \Psi)^\ast \right] \Qiinv G^\ast
\end{split}
\end{equation}
\item Solve, by means of linear algebraic methods
(the Moore-Penrose pseudoinverse), the equation
\begin{equation}
\sum_k \alpha_k Y_k = Y
\end{equation}
\item By linearity, the solution to (\ref{NEWTONSTEP2}) is 
\begin{equation}\label{NEWTONSTEP3}
\Delta\Lambda_i = \sum_k \alpha_k \Sigma_k.
\end{equation}
\end{enumerate}
\end{itemize}

It is clear that the real difficulty here is the computation
of the integrals (\ref{INTEGRAL1}) and (\ref{INTEGRAL2}).
This task requires extensive use of the following results of
 linear stochastic systems theory.

\begin{lemma}
\label{LEMMALYAP}
Let $A$ be a stability matrix and $W(z) = C(zI-A)^{-1}B + D$ 
a minimal realization of a spectral factor of $\Phi(z)$.
Let $\Pi$ be the unique solution to the Lyapunov equation
\begin{equation}
\label{LYAPUNOV}
\Pi = A\Pi A^\ast + BB^\ast
\end{equation}
Then the following hold:
\begin{enumerate}
\item $\int_{-\pi}^{\pi}\Phi(\ejth)\frac{\dth}{2\pi} = C\Pi C^\ast + DD^\ast$.
\item $Z(z) = C(zI-A)^{-1}(A\Pi C^\ast + BD^\ast) + \frac{1}{2}(C\Pi C^\ast + DD^\ast)$
is a realization of the {\em causal part} of $\Phi(z)$; that is,
$Z(z)$ is analytic outside the unit circle and
$\Phi(z) = Z(z) + Z^\ast(z)$.
\end{enumerate}
\end{lemma}

\begin{lemma}
\label{LEMMAARE}
Let $Z(z) = C(zI-A)^{-1}G + \frac{1}{2}\Sigma$ be a
minimal realization of the causal part of a spectrum $\Phi(z)$. Let
$P_-$ be the stabilizing solution of the following
Algebraic Riccati Equation (ARE):
\begin{equation}
\label{ARE}
P = APA^\ast + (G - APC^\ast)(\Sigma - CPC^\ast)^{-1}(G^\ast - CPA^\ast)
\end{equation}
Let moreover $D = (\Sigma - CP_-C^\ast)^{1/2}$ and
$B = (G - AP_-C^\ast)D^{-1}$.
Then $W(z) = C(zI-A)^{-1}B + D$ is the {\em minimum phase}
spectral factor of $\Phi(z)$; that is, $W(z)$ is stable and
with stable causal inverse, and $\Phi(z) = W(z)W^\ast(z)$.
\end{lemma}

\begin{lemma}
\label{LEMMAINVERSE}
Let $F(z) = C(zI-A)^{-1}B + D$ be a square transfer function, where
$D$ is invertible. Then
\begin{equation} 
F^{-1}(z) = -D^{-1}C \left(zI - (A-BD^{-1}C) \right)^{-1}BD^{-1} + D^{-1}
\end{equation}
is a realization of its inverse.
\end{lemma}

\begin{lemma}
\label{LEMMAFERRANTECOLANERI}
For all matrices $P=P^\ast \in \C^{n\times n}$ the following identity holds:
\begin{equation}
\label{EQFC}
\begin{split}
&\left[
\begin{array}{cc}
B^\ast(z^{-1} I - A^\ast)^{-1} & I
\end{array}
\right]
\left[
\begin{array}{cc}
A^\ast P A - P & A^\ast PB \\
B^\ast P A     & B^\ast PB
\end{array}
\right]
\times\\
&\times
\left[
\begin{array}{c}
(zI - A)^{-1}B \\
I
\end{array}
\right] = 0
\end{split}
\end{equation}
\end{lemma}

\begin{lemma}
\label{LEMMARIGHTLEFT}
Let $A$ be a stability matrix and $H(z) = C(zI - A)^{-1}B + D$ be a  minimal realization.
Let $P$ be the solution of the Lyapunov equation
\begin{equation}
\label{LELSDSFS}
P = A^\ast P A + C^\ast C.
\end{equation}
Let $\left[\begin{matrix}K \cr J \end{matrix}\right]$ be  an ortho-normal basis of 
$\ker \left[\begin{matrix} A^\ast P^{1/2} & C^\ast \end{matrix}\right]$ i.e.
\begin{equation}
\label{ONBKER}
\left[\begin{matrix} A^\ast P^{1/2} & C^\ast \end{matrix}\right] 
\left[\begin{matrix}K \cr J \end{matrix}                 \right] = 0, \qquad
\left[\begin{matrix} K^\ast & J^\ast \end{matrix}        \right]
\left[\begin{matrix}K \cr J \end{matrix}                 \right] =I.
\end{equation}
Let $G:=P^{-1/2}K$ and 
\begin{equation}
H_1(z) := (D^\ast C + B^\ast P A) (zI - A)^{-1} G + B^\ast P G + D^\ast J.
\label{KDEF}
\end{equation}
Then, $H^\ast(z) H(z) = H_1(z) H_1^\ast(z)$.
\end{lemma}

Lemmas \ref{LEMMALYAP}, \ref{LEMMAARE} and \ref{LEMMAINVERSE}
are standard results (see for example \cite{FCG}). 
The proofs of Lemmas \ref{LEMMAFERRANTECOLANERI} 
and \ref{LEMMARIGHTLEFT} can be found
in \cite[Appendix A]{FERRANTECOLANERI} and
\cite [Appendix A]{FPRMULTIV}, respectively.

\begin{remark}
\label{REMARKRIGHTLEFT}
Lemma \ref{LEMMARIGHTLEFT} not only gives us a tool to compute a left 
factor from a right  factor of a given spectrum.
It also works in the opposite direction.
Indeed, let $W(z) = C(zI-A)^{-1}B + D$ be a minimal realization,
and let $\zeta = z^{-1}$.
Then 
{\small
\begin{displaymath}
\begin{split}
&\Phi(z) = W(z)W^\ast(z) \\
&= ( C(zI-A)^{-1}B + D ) ( B^\top(z^{-1}I-A^\top)^{-1}C^\top + D^\top ) \\
&= ( B^\top(\zeta^{-1}I-A^\top)^{-1}C^\top + D^\top )^\top ( B^\top(\zeta I-A^\top)^{-1}C^\top + D^\top ) \\
&= ( B^\top(\zeta I-A^\top)^{-1}C^\top + D^\top )^\ast ( B^\top(\zeta I-A^\top)^{-1}C^\top + D^\top ) \\
&:= H^\ast(\zeta) H(\zeta)
\end{split}
\end{displaymath}
}
Applying Lemma \ref{LEMMARIGHTLEFT} we can find an
$H_1(\zeta) = H(\zeta I - F)^{-1}G + K$
such that $H^\ast(\zeta) H(\zeta) = H_1(\zeta) H_1^\ast(\zeta)$.
Now turning back to $z$:
{\small
\begin{displaymath}
\begin{split}
&H_1(\zeta) H_1^\ast(\zeta)\\
&= ( H(\zeta I - F)^{-1}G + K ) ( G^\top(\zeta^{-1} I - F^\top)^{-1}H^\top + K^\top ) \\
&= ( G^\top(z^{-1} I - F^\top)^{-1}H^\top + K^\top )^\top ( G^\top(zI - F^\top)^{-1}H^\top + K^\top ) \\
&= ( G^\top(zI - F^\top)^{-1}H^\top + K^\top )^\ast ( G^\top(zI - F^\top)^{-1}H^\top + K^\top ) \\
&= W_1^\ast(z) W_1(z) = \Phi(z)
\end{split}
\end{displaymath}
}
\end{remark}

\subsection{Factorization of $\Ql^{-1}(z)$}

The first problem to solve is to obtain a spectral factor of
$\Qinv(z)$, where
$\Ql(z) = (I + G^\ast(z)\Lambda G(z))$. To this end,
note that
\begin{equation}
\label{Q1}
\begin{split}
\Ql(z)
&= 
\left[
\begin{array}{cc}
B^\ast(z^{-1} I - A^\ast)^{-1} & I
\end{array}
\right]
\left[
\begin{array}{cc}
\Lambda & 0 \\
0       & I
\end{array}
\right]
\times\\
&\times
\left[
\begin{array}{c}
(zI - A)^{-1}B \\
I
\end{array}
\right]
\end{split}
\end{equation}
Applying lemma \ref{LEMMAFERRANTECOLANERI},
we can rewrite (\ref{Q1}) as
\begin{equation}
\begin{split}
\Ql(z) 
&= 
\left[
\begin{array}{cc}
B^\ast(z^{-1} I - A^\ast)^{-1} & I
\end{array}
\right]
\times \\
&\times
\left[
\begin{array}{cc}
A^\ast P A - P +\Lambda & A^\ast PB \\
B^\ast P A     & B^\ast PB + I
\end{array}
\right]
\times \\
&\times
\left[
\begin{array}{c}
(zI - A)^{-1}B \\
I
\end{array}
\right]
\end{split}
\end{equation}
Now, the following linear matrix inequality:
\begin{equation}
\label{LMI}
\begin{split}
&\left[
\begin{array}{cc}
A^\ast P A - P +\Lambda & A^\ast PB \\
B^\ast P A     & B^\ast PB + I
\end{array}
\right] \\
&=
\left[
\begin{array}{c}
M^\ast  \\
N^\ast 
\end{array}
\right] 
\left[
\begin{array}{cc}
M & N
\end{array}
\right] 
\geq 0
\end{split}
\end{equation}
is solvable for $P = P^\ast > 0$ if and only if such is the following ARE:
\begin{equation}
\label{RICCATI}
P = A^\ast P A - A^\ast PB (B^\ast PB + I)^{-1} B^\ast P A +\Lambda
\end{equation}
The stabilizing solution $P$ of (\ref{RICCATI}) gives a realization
for the square, minimum phase {\em co-analytic} spectral factor of $Q(z)$.
We have:
\newcommand{\Deltal}{{\Delta_\Lambda}}
\begin{equation}
\label{DELTA1}
\begin{split}
N &= N^\ast = (B^\ast PB + I)^{1/2} \\
M &= (B^\ast PB + I)^{-1/2} B^\ast PA \\
\Deltal(z) &=  \left[\begin{array}{cc}M & N\end{array}\right] 
               \left[\begin{array}{c}(zI - A)^{-1}B \\I\end{array}\right] \\
          &= (B^\ast PB + I)^{-1/2} B^\ast PA (zI-A)^{-1} B \\
          &+ (B^\ast PB + I)^{1/2} \\ 
\Ql(z) &= \Deltal^\ast(z) \Deltal(z) 
\end{split}
\end{equation}
and finally $\Ql^{-1}(z) = \Deltal^{-1}(z) \Deltal^{-\ast}(z)$ where,
by means of lemma \ref{LEMMAINVERSE},
\begin{equation}
\label{DELTAINVERSE}
\begin{split}
\Deltal^{-1}(z) &= - (B^\ast PB + I)^{-1} B^\ast PA (zI-\Gamma)^{-1} B \times \\
               &\ \ \times(B^\ast PB + I)^{-1/2} + (B^\ast PB + I)^{-1/2} \\
\Gamma &= A - B (B^\ast PB + I)^{-1} B^\ast PA
\end{split}
\end{equation}

\subsection{Computation of the integrals in (\ref{INTEGRAL1}) and (\ref{INTEGRAL2})}

\newcommand{\Phione}{{\Phi_\Psi}}

By virtue of Lemma \ref{LEMMARIGHTLEFT}
and Remark \ref{REMARKRIGHTLEFT}, we
can switch from a {\em right} factorization of a spectrum
($\Phi=H^\ast H$) to a {\em left} factorization ($\Phi=WW^\ast$),
and vice versa.
We will now show that both (\ref{INTEGRAL1}) and (\ref{INTEGRAL2}) can 
be reduced to integrals of the form
\begin{equation}
\label{INTEGRALFORM}
\int G(z) \Deltal^{-1}(z) \Phi(z) \Deltal^{-\ast}(z) G^\ast(z)
\end{equation}
where $\Phi(z)$ is a spectrum.
Indeed, let $\Phione(z) = \Deltal^{-\ast}(z)\Psi(z)\Deltal^{-1}(z)$.
Then
\begin{equation}
\label{INTEGRALFORM2}
\begin{split}
&\int G(z) \Qinv(z) \Psi(z) \Qinv(z) G^\ast(z) \\
&= \int G(z) \Deltal^{-1}(z) \left( \Deltal^{-\ast}(z) \Psi(z) \Deltal^{-1}(z) \right) \Deltal^{-\ast}(z) G^\ast(z)
\end{split}
\end{equation}
which has the form (\ref{INTEGRALFORM})
with $\Phi = \Phione$.
Applying Lemma \ref{LEMMARIGHTLEFT}
we obtain a (left) spectral factor of $\Phione(z)$:
\begin{equation}
\begin{split}
\Phione(z) 
&= \Deltal^{-\ast}(z) (W_\Psi(z) W_\Psi^\ast(z)) \Deltal^{-1}(z) \\
&= \Deltal^{-\ast}(z) (H^\ast_\Psi(z) H_\Psi(z)) \Deltal^{-1}(z) \\
&= (H_\Psi(z)\Deltal^{-1}(z))^\ast (H_\Psi(z)\Deltal^{-1}(z)) \\
&= W_1(z) W_1^\ast(z)
\end{split}
\end{equation}
Finally,
(\ref{INTEGRALFORM2}) can be computed obtaining a realization
of $G(z) \Deltal^{-1}(z) W_1(z)$ and applying Lemma \ref{LEMMALYAP}.
Now, let $\Phi_\Sigma(z) = \Deltal^{-\ast}(z)G^\ast(z) \Sigma G(z)\Deltal^{-1}(z)$,
where $\Sigma$ is one of the precomputed generators of $\Rgamma$.
Then
%
%
%
%
%
%
\begin{equation}
\label{INTEGRAL11}
\begin{split}
&\int G \Qinv \left[ (G^\ast \Sigma G \Qinv \Psi) + (G^\ast \Sigma G \Qinv \Psi)^\ast \right] \Qinv G^\ast \\
&= \int G \Deltal^{-1}\Deltal^{-\ast} 
   \left[ (G^\ast \Sigma G \Deltal^{-1}\Deltal^{-\ast} \Psi) \right. \\
   &\quad\quad\quad\quad \left.+\ (\Psi \Deltal^{-1}\Deltal^{-\ast}  G^\ast \Sigma G ) \right] \Deltal^{-1}\Deltal^{-\ast} G^\ast \\
&= \int G \Deltal^{-1} 
   \left[ (\Deltal^{-\ast}G^\ast \Sigma G \Deltal^{-1}\Deltal^{-\ast} \Psi\Deltal^{-1}) \right.\\
   &\quad\quad\quad\quad \left.+\ (\Deltal^{-\ast}\Psi \Deltal^{-1}\Deltal^{-\ast}  G^\ast \Sigma G \Deltal^{-1}) \right] \Deltal^{-\ast} G^\ast \\ 
&= \int G \Deltal^{-1} \left[ \Phi_\Sigma \Phione + \Phione \Phi_\Sigma \right] \Deltal^{-\ast} G^\ast  \\
&= \int G \Deltal^{-1} \left[ (\Phi_\Sigma + \Phione) (\Phi_\Sigma + \Phione) \right.\\
   &\quad\quad\quad\quad \left.-\ \Phione \Phione - \Phi_\Sigma \Phi_\Sigma \right] \Deltal^{-\ast} G^\ast  \\
&= \int G \Deltal^{-1} \left[ (\Phi_\Sigma + \Phione) (\Phi_\Sigma + \Phione) \right] \Deltal^{-\ast} G^\ast \\
&- \int G \Deltal^{-1} \left[ \Phione \Phione \right] \Deltal^{-\ast} G^\ast 
 - \int G \Deltal^{-1} \left[ \Phi_\Sigma \Phi_\Sigma \right] \Deltal^{-\ast} G^\ast 
\end{split}
\end{equation}
which is a difference of integrals of the form (\ref{INTEGRALFORM}).
To compute (\ref{INTEGRAL11}), we must obtain (left) spectral
factors of $\Phione \Phione^\ast$, $\Phi_\Sigma \Phi_\Sigma^\ast$ and
$(\Phi_\Sigma + \Phione)(\Phi_\Sigma + \Phione)^\ast$.
Suppose, first, that $\Sigma>0$.
For the first spectrum we have
\begin{equation}
\begin{split}
\Phione \Phione^\ast 
&= W_1 (W_1^\ast W_1) W_1^\ast
= W_1 (H_1 H_1^\ast) W_1^\ast\\
&= (W_1 H_1) (W_1 H_1)^\ast
\end{split}
\end{equation}
For the second, we have 
\begin{equation}
\begin{split}
\Phi_\Sigma  
&= (\Deltal^{-\ast}G^\ast \Sigma^{1/2}) (\Sigma^{1/2} G\Deltal^{-1})
= H_\Sigma^\ast H_\Sigma
= W_\Sigma W_\Sigma^\ast \\
\Phi_\Sigma \Phi_\Sigma^\ast 
&= W_\Sigma (W_\Sigma^\ast W_\Sigma) W_\Sigma^\ast
= W_\Sigma (K_\Sigma K_\Sigma^\ast) W_\Sigma^\ast\\
&= (W_\Sigma K_\Sigma) (W_\Sigma K_\Sigma)^\ast
\end{split}
\end{equation}
And for the third:
\begin{equation}
\begin{split}
&(\Phi_\Sigma + \Phione)(\Phi_\Sigma + \Phione)^\ast\\
&= (Z_\Sigma + Z_\Sigma ^\ast + Z_1 + Z_1^\ast)(Z_\Sigma + Z_\Sigma ^\ast + Z_1 + Z_1^\ast)^\ast \\
&= ((Z_\Sigma + Z_1) + (Z_\Sigma + Z_1)^\ast)((Z_\Sigma + Z_1) + (Z_\Sigma + Z_1)^\ast)^\ast \\
&= (Z_{1\Sigma} + Z_{1\Sigma}^\ast)(Z_{1\Sigma} + Z_{1\Sigma}^\ast)^\ast \\
&= W_{1\Sigma} (W_{1\Sigma}^\ast W_{1\Sigma}) W_{1\Sigma}^\ast \\
&= W_{1\Sigma} (H_{1\Sigma} H_{1\Sigma}^\ast) W_{1\Sigma}^\ast \\
&= (W_{1\Sigma} H_{1\Sigma}) (W_{1\Sigma} H_{1\Sigma})^\ast 
\end{split}
\end{equation}
where $Z_1$ is the causal part of $\Phione$,
$Z_{1\Sigma} = Z_1 + Z_\Sigma$, 
$W_{1\Sigma}$ 
is a left factor of the spectrum
$Z_{1\Sigma} + Z_{1\Sigma}^\ast$,
$H_{1\Sigma}$ 
is a left factor of the spectrum
$W_{1\Sigma}^\ast W_{1\Sigma}$,
and
where we used Lemma \ref{LEMMALYAP} to obtain the causal part
of $\Phi_\Sigma$ and $\Phione$ from their spectral factors,
and Lemma \ref{LEMMAARE} to obtain the minimum phase spectral
factor of the sum $\Phi_\Sigma + \Phione$ from its causal part.
Thus, if $\Sigma>0$, we really have all the tools to compute integral
(\ref{INTEGRAL11}).
%
%

Now, $\Sigma$ is {\em not} necessarily positive definite, but
if $-\lambda < 0$ is the minimum between the eigenvalues of all the
generators $\Sigma_k$, then $\Sigma + (\lambda + 1) I$ {\em is} positive
definite. Thus, in the general case, by linearity
(\ref{INTEGRAL11}) can be reduced to:
\begin{equation}
\begin{split}
&\int G \Deltal^{-1} \left[ \Phi_\Sigma \Phione + \Phione \Phi_\Sigma \right] \Deltal^{-\ast} G^\ast  \\
&=\int G \Deltal^{-1} \left[ \Phi_{\Sigma + (\lambda+1)I - (\lambda+1)I}\ \Phione \right. \\
&\quad\quad\quad\quad\quad     \left.+\ \Phione\ \Phi_{\Sigma + (\lambda+1)I - (\lambda+1)I} \right] \Deltal^{-\ast} G^\ast  \\
&=\int G \Deltal^{-1} \left[ \Phi_{\Sigma + (\lambda+1)I}\ \Phione + \Phione\ \Phi_{\Sigma + (\lambda+1)I} \right] \Deltal^{-\ast} G^\ast  \\
&-(\lambda+1)\int G \Deltal^{-1} \left[ \Phi_I\ \Phione + \Phione\ \Phi_I \right] \Deltal^{-\ast} G^\ast  \\
\end{split}
\end{equation}
which is a difference of integrals 
with the same structure of (\ref{INTEGRAL11}), and
that are computable with the above tools 
(obviously 
$\int G \Deltal^{-1} \left[ \Phi_I\ \Phione + \Phione\ \Phi_I \right] \Deltal^{-\ast} G^\ast$
needs to be computed only once).
This enables us to solve equation (\ref{NEWTONSTEP}).

\subsection{Computations in the backtracking step}

The backtracking stage involves similar, though easier, computations.
We must check the following conditions:
\begin{equation}
\label{BACKTRACKINGDOMAINJ}
\Lambda_{i} + t_i^k\Delta \Lambda_i \in \LHG
\end{equation}
\begin{equation}
\label{BACKTRACKINGJ}
\Jpsi(\Lambda_i + t_i^k\Delta \Lambda_i) < \Jpsi(\Lambda_i) + \alpha t_i^k \nabla \Jpsi_{\Lambda_i} \Delta \Lambda_i
\end{equation}

Checking (\ref{BACKTRACKINGDOMAINJ}) is really a matter of checking
whether we can factorize $I + G^\ast (\Lambda_i + t_i^k \Delta\Lambda_i) G$.
Thus $t_i^k$ must be halved until the ARE
(\ref{RICCATI}) is solvable having $\Lambda = \Lambda_i + t_i^k \Delta\Lambda_i$.

Finally, to check (\ref{BACKTRACKINGJ}), we need to compute $\Jpsi$.
This can be done in a way similar to the above computations:
\begin{equation}
\label{JPSICOMPUTATION}
\begin{split}
\Jpsi(\Lambda) 
&= \tr \int (I + G^\ast \Lambda G)^{-1} \Psi + \tr \Lambda \\
&= \tr \int \Deltal^{-1}\Deltal^{-\ast} W_\Psi W_\Psi^\ast + \tr \Lambda \\
&= \tr \int \Deltal^{-\ast} (W_\Psi W_\Psi^\ast) \Deltal^{-1} + \tr \Lambda \\
&= \tr \int \Deltal^{-\ast} (H_\Psi^\ast H_\Psi) \Deltal^{-1} + \tr \Lambda \\
&= \tr \int (H_\Psi \Deltal^{-1})^\ast (H_\Psi \Deltal^{-1}) + \tr \Lambda \\
&= \tr \int W W^\ast + \tr \Lambda
\end{split}
\end{equation}
%

\section{Proof of global convergence}
\label{SECTIONCONVERGENCE}

Given that the minimum of $\Jpsi$ exists {\em and is unique}, 
we investigate global convergence of 
our Newton algorithm. First, we recall the following

{\em Definition}:
a function $f(x)$ twice differentiable in a set $S$ is said to be
{\em strongly convex} in $S$ if there exists a constant $m>0$
such that $H(x) \geq mI$ for $x \in S$, where $H(x)$ is the Hessian of $f$ at $x$.\\

\noindent
We restrict our analysis to a sublevel set of $\Jpsi$.
Let $\Lambda_0 = 0$. The set
\begin{equation}
S := \left\{\Lambda \in \LHG : \Jpsi(\Lambda) \leq \Jpsi(\Lambda_0) = \tr \smallint \Psi\right\}
\end{equation}
is {\em compact} (as it was shown in \cite[Section VII]{FPRMULTIV}). Because of the backtracking in the algorithm,
the sequence $\Jpsi(\Lambda_0), \Jpsi(\Lambda_1), ...$ is  decreasing. Thus $\Lambda_n\in S,\,\forall n\ge 0$. We now wish to apply a theorem in \cite[9.5.3, p. 488]{BV} on convergence of the Newton algorithm with backtraking for strongly convex functions on $\R^n$.
This theorem ensures linear decrease for a finite
number of steps, and quadratic convergence to the minimum after the linear stage,
thus establishing {\em global} convergence of the Newton algorithm
with backtracking.
We proceed to establish first {\em strong convexity} of $J_\Psi$ on $S$.
To do that, we employ the following result.
\begin{lemma}
\label{STRONGCONVEXITY}
Let $f(x)$ be defined over an open convex subset $D$ of a finite-dimensional
linear space $V$. Assume that $f$ is twice continuously differentiable
and strictly convex on $D$. Then $f$ is {\em strongly}
convex on any compact set $S \subset D$.
\end{lemma}

\IEEEproof
 
First, recall that since $f$ is twice continuously differentiable
and strictly convex, its Hessian $H_x$ is an Hermitian positive-definite
matrix at each point $x$. By Lemma \ref{LEMMA0}, the mapping from $H$ to its minimum
(real) eigenvalue is continuous. 
It follows that the mapping from $x$ to the minimum eigenvalue of
the Hessian of $f$ at $x$ is also continuous, 
being a composition of continuous functions. Hence the latter admits a minimum $m$ in the compact set $S$ by Weierstrass' theorem.
Thus $m$ is the minimum of the eigenvalues of all the Hessians computed
in $S$, and $m$ cannot be zero, since otherwise there would be an $x$
with $H_x$ singular, and this cannot happen since $f$ is strictly convex.
Hence $H_x \geq mI \ \forall x \in S$, i.e. $f$ is strongly
convex on $S$.
\qed

\begin{remark}
By an argument similar to that of Lemma \ref{LEMMA0},
it can be shown that for a twice continuously differentiable
function which is strictly convex on $D$, 
there exists $M>0$ such that $H_x \leq MI$ for all $x\in S$.
Moreover, strong convexity on a {\em closed} set $S$ implies
boundedness of the latter.
Thus, strong convexity and boundedness of the Hessian are
intertwined, and both are {\em essential} in the proof of 
Theorem \ref{BVTHEOREM} (see \cite{BV}).
\end{remark}

\begin{theorem}
\label{BVTHEOREM}
The following facts hold true:
\begin{enumerate}
\item $\Jpsi$ is twice {\em continuously} differentiable on $S$;
\item $\Jpsi$ is strongly convex on $S$;
\item the Hessian of $\Jpsi$ is Lipschitz-continuous over $S$;
\item the sequence $\{\Lambda_i; i\ge 0\}$ generated by the Newton algorithm of Section V (\ref{NEWTONSTEP0})-(\ref{BACKTRACKINGJ}) converges to the unique minimum point of $J_\Psi$ in $\LHG$.
\end{enumerate} 
\end{theorem}
\IEEEproof

Property 1 is a trivial consequence of Theorem \ref{JPSIREGULARITYPROP}.
To prove  2, remember that $\Jpsi$ is strictly convex on $\LHG$, hence
also on $S$, and apply Lemma \ref{STRONGCONVEXITY}.
As for property 3, what it really says is that the following operator:
\begin{displaymath}
{\cal H}: \Lambda \mapsto \Hess(\cdot, \cdot)
\end{displaymath}
is Lipschitz continuous on $S$.
Theorem \ref{JPSIREGULARITYPROP} implies that $\Jpsi \in C^3(\LHG)$ or,
which is the same, that ${\cal H} \in C^1(\LHG)$. The continuous
differentiability of ${\cal H}$ implies its Lipschitz continuity
over an arbitrary compact subset of $\LHG$, hence also over
the sublevel set $S$, and property 3 follows. 
Finally, to prove 4, notice that all the hypotheses of \cite[9.5.3, p. 488]{BV} are satisfied. Namely, the function to be minimized $\Jpsi$ is strongly convex on the compact set $S$, and its Hessian is  Lipschitz-continuous over $S$. It remains to observe that $\Jpsi$ is defined
over a subset of the linear space $\Rgamma$ which has {\em finite dimension} $d$ over $\R$
(recall that $\Rgamma$ is spanned by a finite set of matrices.
See Proposition \ref{RANGEGAMMAPROP} and Remark \ref{RANGEGAMMAREMARK},
where $d \leq N$).
Thus, once we choose a base in $\Rgamma$, to every $\Lambda \in \LHG$ there corresponds a vector in $\R^d$, to every positive definite bilinear form over $\Rgamma$ there corresponds
a positive definite matrix in $\R^{d \times d}$, and to every compact set in $\LHG$ there corresponds a compact set in $\R^d$. Hence, every convergence result that holds in $\R^d$ must also hold in the abstract setting, in view of the homeomorphism between one space and the other.
\qed

\section{Application to spectrum estimation}

\subsection{A spectral estimation procedure}\label{sep}

Following the purposes of the THREE method presented in \cite{BGLTHREE},
now we describe an application of the above approximation algorithm to the
estimation of spectral densities.
Consider  first the scalar case, and suppose that the
finite sequence $y_1, ..., y_N$ is
extracted from a realization of a zero-mean, weakly stationary discrete-time 
process $\{y_t\}_{t=-\infty}^{+\infty}$.
We want to estimate the spectral density $\Phi_y(\ejth)$ of $y$.
The idea is the following:

\begin{itemize}
\item Fix a transfer function $G(z) = (zI-A)^{-1}B$, feed the data $\{y_i\}$ to it,
and collect the output data $\{x_i\}$.
\item Compute a consistent, and possibly unbiased, estimate $\hat{\Sigma}$ of
the covariance matrix of the outputs $\{x_i\}$. Note that some output samples 
$x_1, ..., x_M$ should be discarded so that the filter can be considered to operate in steady state.
\item Choose as  ``prior'' spectrum $\hat{\Phi}_y$ a coarse, low-order, estimate
of the true spectrum of $y$ obtained by means of another (simple) identification method.
\item ``Refine'' the estimate $\hat{\Phi}_y$ by
solving the approximation problem (\ref{PRIMALPROBLEM}) with respect
to $G(z)$, $\Sigma = \hat{\Sigma}$, and $\Psi = \hat{\Phi}_y$.
\end{itemize}

To be clear, the result of the above procedure is {\em the only spectrum, compatible
with the output variance $\hat{\Sigma}$, which is closest to the rough estimate
$\hat{\Phi}_y$ in the $d_H$ distance}.
Note that we are left with significant degrees of freedom in applying the above procedure:
The method for estimating $\hat{\Phi}_y$, in particular its degree, and
the whole structure of $G(z)= (zI-A)^{-1}B$, which has no contraints other than  $A$ being
a stability matrix and  $(A,B)$ being reachable.

The coarsest possible estimate of $\Phi_y$ is the 
constant spectrum equal to the sample variance of the $\{y_i\}$, i.e.\ 
$\hat{\Phi}_y(\ejth) \equiv \hat{\sigma}_y^2$, where
$\hat{\sigma}_y^2 = \frac{1}{N-1}\sum_{i=1}^N |y_i|^2$.
The resulting spectrum has the form 
$\hat{\sigma}_y^2 (1 + G^\ast \hat{\Lambda} G)^{-2}$.
Another simple choice is $\hat{\Phi}_y = W(z)W^\ast(z)$,
where $W(z) = \hat{\sigma}_e \frac{c(z)}{a(z)}$
is a low-order AR, MA or ARMA model estimated from $y_1, ..., y_N$ by means of 
predictive error minimization methods or the like.

The flexibility in the choice of $G(z)$ is more essential, and has
more profound implications. As described in \cite{BGLTHREE}, \cite{GLKL},
\cite{FPRHELL} and \cite{FPRMULTIV}, the following choice:
\newcommand{\bmat}{\left[ \begin{array}}
\newcommand{\emat}{\end{array} \right]}
\begin{equation}
A=\bmat{ccccc}
p_1&0&0&\dots&0\\
0&p_2&0&\dots&0\\
\vdots&\vdots&
&\ddots&\vdots\\
       0&0&0&\dots&0\\
       0&0&0&\dots&p_n
       \emat,\quad
B=\bmat{c}1\\1\\\vdots\\1\\1\emat
\end{equation}
where the $p_i$'s lie inside the unit circle,
implies that the (true) steady-state variance $\Sigma$ has the structure of 
a Pick matrix, and the corresponding problem of finding {\em any} spectrum that satisfies
(\ref{CONSTRAINT}) is a Nevanlinna-Pick interpolation. 
Moreover, the following choice:
\begin{equation}
\label{COVARIANCEEXTENSION}
A = \bmat{ccccc}
     0&1&0&\dots&0\\
     0&0&1&\dots&0\\
     \vdots&\vdots&
     &\ddots&\vdots\\
     0&0&0&\dots&1\\
     0&0&0&\dots&0
     \emat, \quad 
B =\bmat{c}0\\0\\\vdots\\0\\1\emat
\end{equation}
implies that the steady-state variance $\Sigma$ is a Toeplitz matrix
whose diagonals contain the lags $c_0, c_1, ..., c_{n-1}$ of the covariance
signal of the input,
and the corresponding problem of finding {\em any} spectrum that satisfies
(\ref{CONSTRAINT}) is a covariance extension problem. 

These facts justify the theoretical interest in algorithms for constrained
spectrum approximation, if for no other reason,
as tools to compute at least {\em one} solution to a Nevanlinna-Pick interpolation
or to a covariance extension problem, respectively.
But the freedom in choosing $G(z)$ has implications also in the above
practical application to spectral estimation, where the key properties,
not surprisingly, depend on the poles of $G(z)$, i.e., the eigenvalues of $A$.
In general, as described in \cite{BGLTHREE},
the {\em magnitude} of the latter has implications
on the variance of the sample covariance $\hat{\Sigma}$: The closer
the eigenvalues to the origin, the smaller that variance 
(see \cite[Section II.D]{BGLTHREE}).
Moreover, at least as far as THREE \cite{BGLTHREE} is concerned, the {\em phase} of the
eigenvalues influences resolution capability:
More precisely, {\em the spectrum estimation procedure has higher resolution
in those sectors of the unit circle where more eigenvalues are located}.
According to simulations, the latter statement appears to be true
also in our setting (the fundamental difference being that the metric which
is minimized is the Hellinger distance instead of the Kullback-Leibler one).

\begin{remark}
In the above setting $\hat{\Sigma}$ is 
a consistent estimate of the true steady-state variance.
Although $\hat{\Sigma}$ must belong to $\Rgamma$
as $N\rightarrow +\infty$
(this being the case even if $y$ is the sum of a purely nondeterministic
process and some sinusoids, as in the simulations that follow),
it is almost certainly not the case that $\hat{\Sigma} \in \Rgamma$
when we have available only the finitely many data $x_{M+1}, ..., x_N$.
Strictly speaking, this implies that the contraint (\ref{CONSTRAINT})
with $\Sigma = \hat{\Sigma}$ is almost always not feasible.
It turns out that, increasing  the tolerance threshold in its step 5, 
the Newton algorithm exhibits some kind of robustness in this respect.
That is, it leads to a $\Lambda$ whose corresponding spectrum
$\hat{\Phi}$ is  {\em close} to satisfying the constraint.

Nevertheless, we prefer a clear understanding of what the resulting spectrum really is.
Thus, we choose to enforce feasibility of the approximation problem,
at least as permitted by machine number representation,
before starting the optimization procedure.
To this end,  following the same approach employed in \cite{BGLTHREE}, we pose the approximation
problem not in terms of the estimated $\hat{\Sigma}$,
but in terms of its {\em orthogonal projection} $\hat{\Sigma}_\Gamma$ 
onto $\Rgamma$, which can be easily computed  by means of algebraic methods.
That is to say: We cannot approximate in the preimage $\Gamma^{-1}(\hat{\Sigma})$,
because that set is empty, thus we choose to approximate in $\Gamma^{-1}(\hat{\Sigma}_\Gamma)$,
where $\hat{\Sigma}_\Gamma$ is the 
matrix closest to $\hat{\Sigma}$ such that its preimage is not empty.
This seems a reasonable choice and by the way it is,
{\em mutatis mutandis}, what the Moore-Penrose pseudoinverse does
for the ``solution'' $\hat{x} = A^\dagger b$, when the linear system $Ax = b$ is
not solvable.

Note that it is not guaranteed at all that the projection of a positive definite
matrix onto a subspace of the Hermitian matrices is itself positive definite.
In practice, this is not really a problem, inasmuch $\hat{\Sigma}$ is
``sufficiently positive'' and close to $\Rgamma$.
The positivity of $\hat{\Sigma}_\Gamma$ must anyway be checked before proceeding.
This approach and the considerations on the
positivity issue should be compared to \cite[Section II.D]{BGLTHREE}, which
deals with the particular case when $\Rgamma$ is the space of Toeplitz matrices,
and to \cite[Section 4]{GDISTANCESBETWEEN}, where,
to find a matrix a $\hat{\Sigma}_\Gamma$ close to $\hat{\Sigma}$,
a Kullback-Leibler criterion is adopted instead of least squares.
\end{remark}

\subsection{Simulation results: Scalar case}
\newcommand{\ee}{{\rm e}}
\newcommand{\jj}{{\rm j}}
%
\begin{figure}[h!]
\begin{center}
\includegraphics[width=13cm]{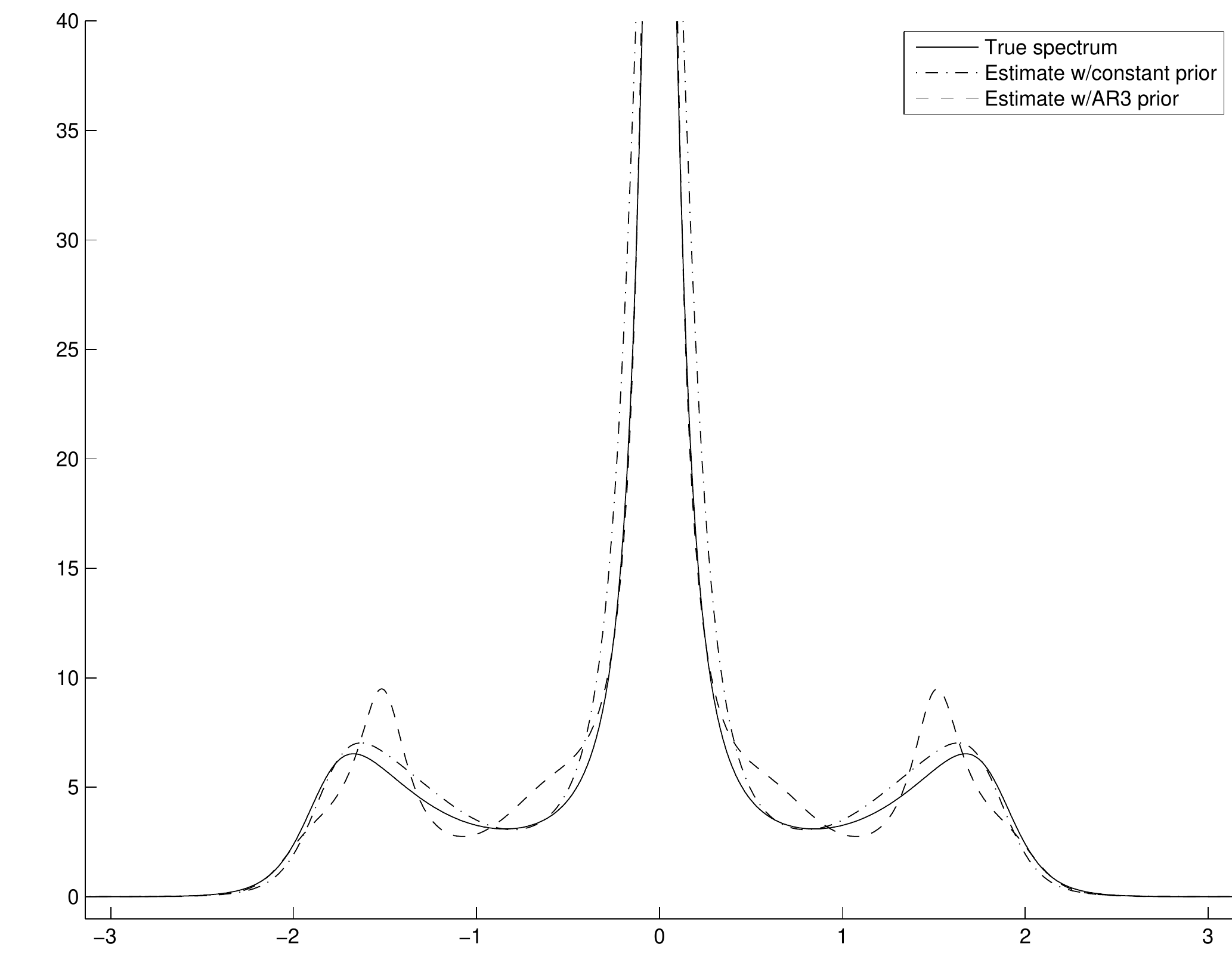}
\end{center}
\caption{Estimation of an ARMA(6,4) spectrum
by means of Hellinger-distance spectrum approximation, constant prior
and AR(3) prior.} 
\label{COVEXTFIGURE}
\end{figure}

Figure \ref{COVEXTFIGURE} shows the results of the above estimation procedure with
$G(z)$ structured according to the covariance extension setting
(\ref{COVARIANCEEXTENSION}) with $6$ covariance lags (i.e.\ $n=6$, $A$ is $6\times 6$),
run over $500$ samples of the following ARMA process:
\begin{displaymath}
\begin{split}
 y(t) &= 0.5 y(t-1) - 0.42 y(t-2) + 0.602 y(t-3) \\
 &- 0.0425 y(t-4) + 0.1192 y(t-5) \\
 &+ e(t) + 1.1 e(t-1) + 0.08 e(t-2) -  0.15 e(t-3) 
\end{split}
\end{displaymath}
(poles in $0.9, -0.2 \pm 0.7\jj, \pm 0.5\jj$)
where $e(t)$ is a zero-mean Gaussian white noise with
unit variance.
Two priors, both estimated from data, have been considered: the constant spectrum
$\hat{\Phi}_y(\ejth) \equiv \hat{\sigma}_y^2$ and the
spectrum $\hat{\Phi}_y = W_{AR}(z)W_{AR}^\ast(z)$, where 
$W_{AR}(z) = \frac{\hat{\sigma}_e}{a(z)}$ is an AR model of order $3$
obtained from the data by means of the Predictive Error Method procedure 
in Matlab's System Identification toolbox.

\begin{figure}[h!]
\begin{center}
\includegraphics[width=13cm]{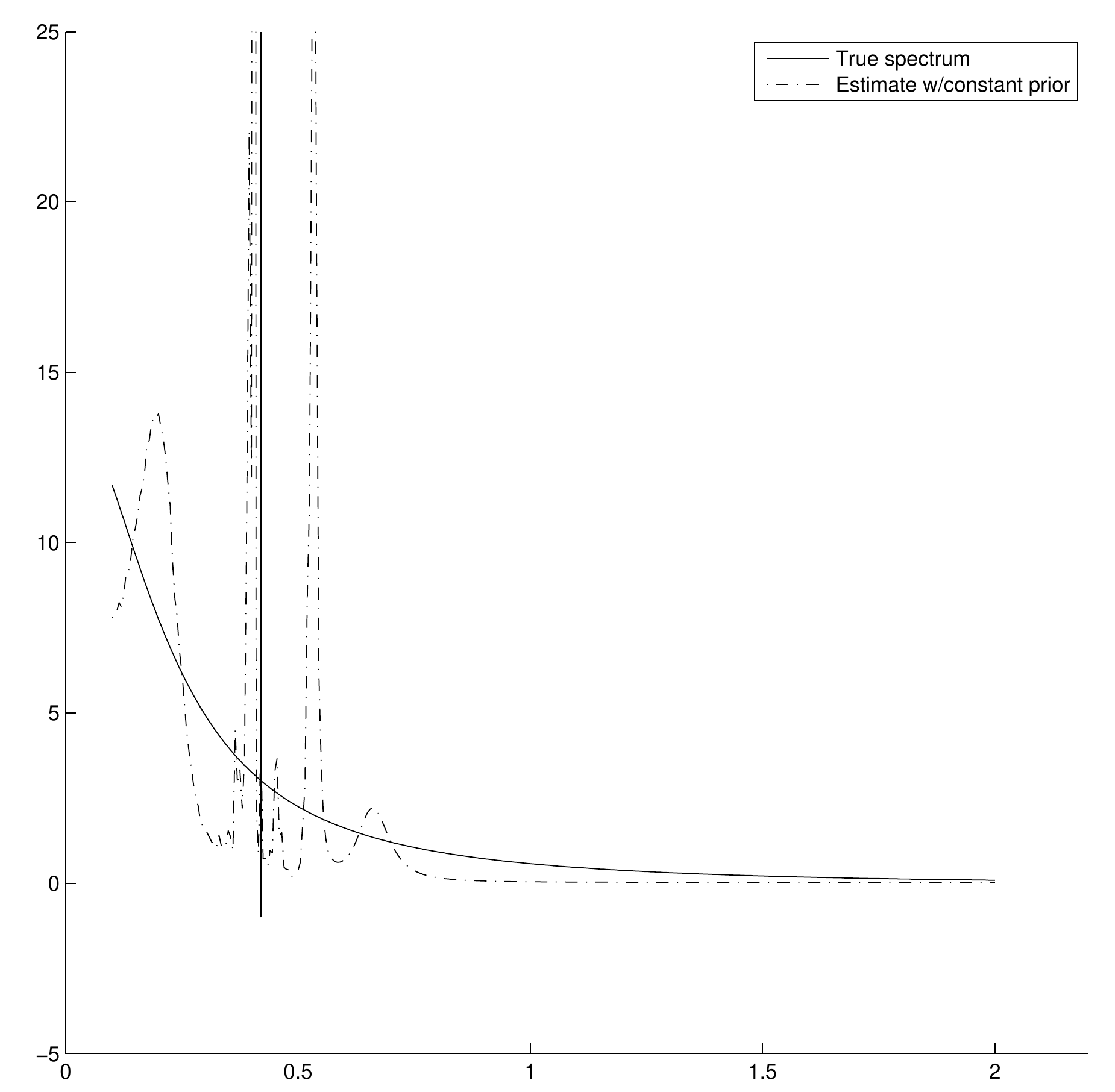}
\end{center}
\caption{Spectral estimates of two sinusoids with superimposed noise
by means of Hellinger-distance spectrum approximation, constant prior.
Compare with \cite[Section IV.B, Example 1]{BGLTHREE}.} 
\label{BGL1FIGURE}
\end{figure}

Figure \ref{BGL1FIGURE} shows the performance of the above procedure in a setting
that resembles that of \cite[Section IV.B, Example 1]{BGLTHREE}.
The estimation procedure was run on 300 samples of a
superposition of two sinusoids in colored noise:
\begin{displaymath}
\begin{split}
y(t) &= 0.5 \sin(\omega_1 t + \phi_1) + 0.5 \sin(\omega_2 t + \phi_2) + z(t) \\
z(t) &= 0.8 z(t-1) + 0.5 \nu(t) + 0.25 \nu(t-1) 
\end{split}
\end{displaymath}
with $\phi_1$, $\phi_2$ and $\nu(t)$ independent normal random variables with zero mean
and unit variance, $\omega_1 = 0.42$ and $\omega_2 = 0.53$.
The prior here considered is the constant spectrum equal to the sample
variance of the $\{y_i\}$ data.
Following \cite{BGLTHREE}, $A$ was chosen real block-diagonal with the following poles
(equispaced in a narrow range where the frequencies of the two sinusoids lie,
to increase resolution in that region):
\begin{displaymath}
\begin{split}
&0, 0.85, -0.85,\\
&0.9 \ee^{\pm\jj 0.42},
0.9 \ee^{\pm\jj 0.44},
0.9 \ee^{\pm\jj 0.46},
0.9 \ee^{\pm\jj 0.48},
0.9 \ee^{\pm\jj 0.50}
\end{split}
\end{displaymath}
(and $B$ a column of ones).
It can be seen that Hellinger-distance based approximation 
does a good job, as does the THREE algorithm, at detecting the spectral
lines at frequencies $\omega_1$ and $\omega_2$.
\subsection{Simulation results: Multivariate case}
\begin{figure}[h!]
\begin{center}
\includegraphics[width=13cm]{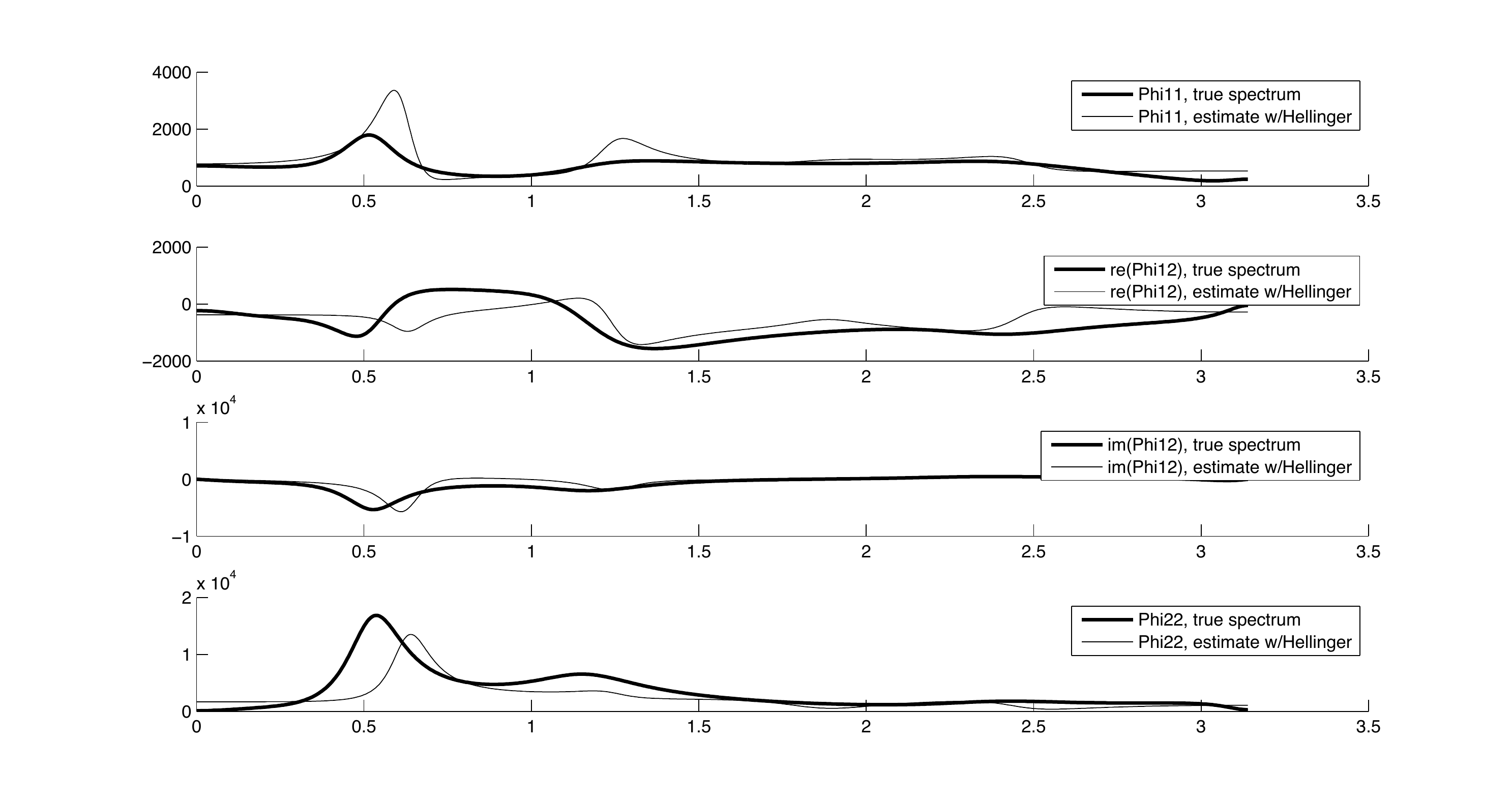}
\end{center}
\caption{Estimation of the spectrum of a bivariate process with rich dynamics
by means of Hellinger-distance spectrum approximation, constant prior.} 
\label{ESTIMFIGURE}
\end{figure}

\begin{figure}[h!]
\begin{center}
\includegraphics[width=13cm]{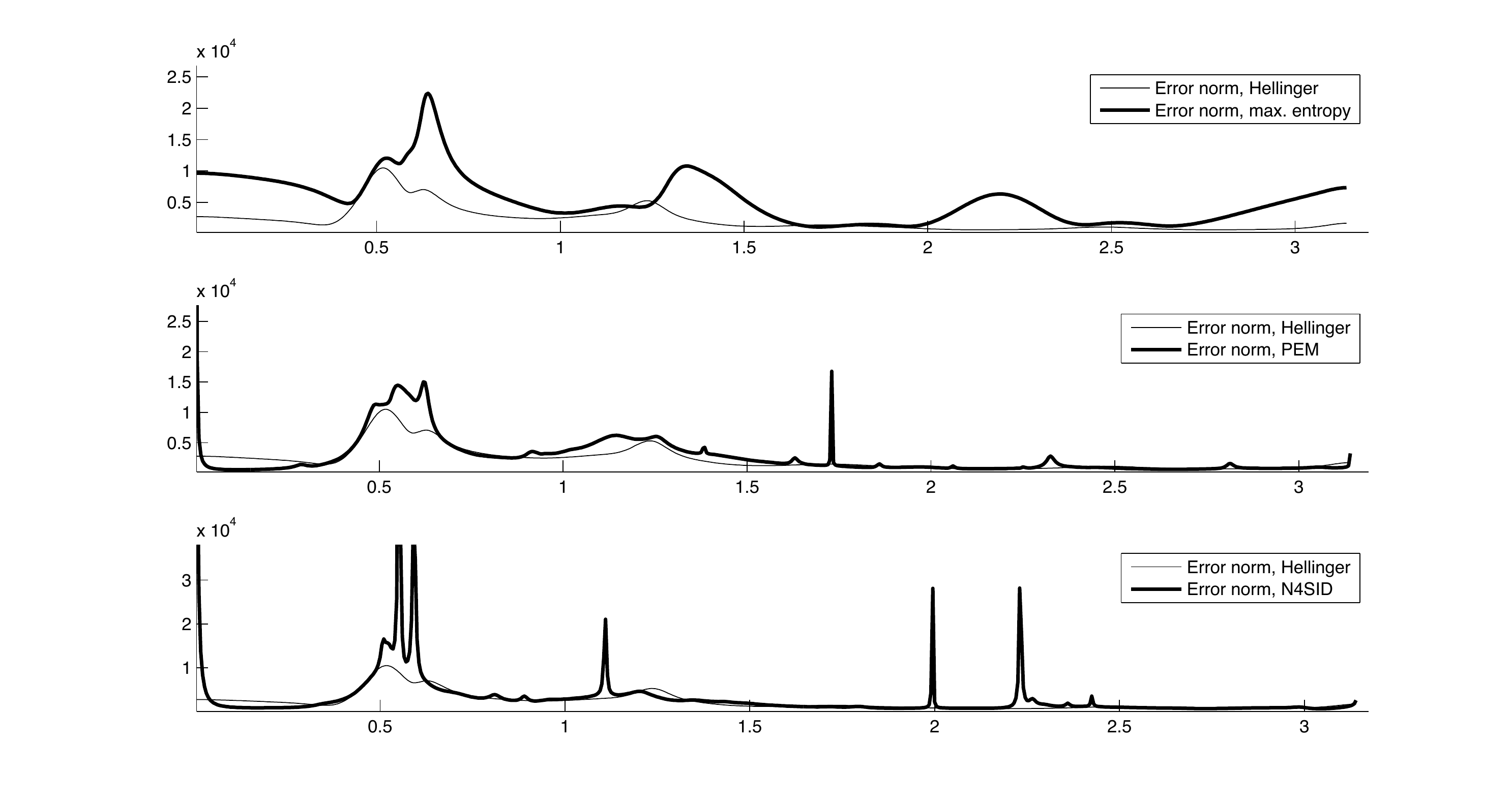}
\end{center} 
\caption{Estimation of the spectrum of a bivariate process with rich dynamics
by means of various methods. Comparison between the spectral norm of the differences
$\hat{\Phi}_{\rm H} - \Phi$,
$\hat{\Phi}_{\rm ME} - \Phi$,
$\hat{\Phi}_{\rm PEM} - \Phi$, and
$\hat{\Phi}_{\rm N4SID} - \Phi$
(average over 50 simulations).}
\label{COMPARISONFIGURE}
\end{figure}

We now consider spectral estimation for a multivariate process. Here, 100 samples
of a bivariate process with a high order spectrum were generated by feeding
a bivariate Gaussian white noise with mean 0 and variance $I$
to a square (stable) shaping filter of order $40$. The latter
was constructed with random coefficients, except for 
one fixed conjugate pair of poles with radius $0.9$ and argument $0.52$, and 
one fixed conjugate pair of zeros with radius $1 - 10^{-5}$ and argument $0.2$.
The transfer function $G(z)$ was chosen with one pole in the origin and
$4$ complex pole pairs with radius $0.9$ and frequencies equispaced in the range $[0, \pi]$.
Then the above estimating procedure was applied, 
with prior spectrum 
chosen as the constant density 
equal to the sample covariance of the bivariate process $y$.
Figure \ref{ESTIMFIGURE} shows a plot of $\Phi_{11}(\ejth)$,
$\Real \Phi_{12}(\ejth)$, $\Imag \Phi_{12}(\ejth)$ and $\Phi_{22}(\ejth)$,
respectively for the true spectrum and for the estimation of the latter
based on one run of 100 samples.
%
%
In Figure \ref{COMPARISONFIGURE}
we compare the performances of various spectral estimation methods in the following way.
We consider four estimates $\hat{\Phi}_{\rm H}$, $\hat{\Phi}_{\rm ME}$,
$\hat{\Phi}_{\rm PEM}$, and $\hat{\Phi}_{ \rm N4SID}$ of $\Phi$.
The spectral density $\hat{\Phi}_{\rm H}$ is the estimate obtained by  the procedure  described above in Subsection \ref{sep}.
The spectral density $\hat{\Phi}_{\rm ME}$ is the maximum entropy estimate \cite{G5} obtained using the same $G(z)$ employed to obtain our estimate.
The spectral densities $\hat{\Phi}_{\rm PEM}$ and $\hat{\Phi}_{\rm N4SID}$ are the estimates of $\Phi$
obtained by using ``off-the-shelf'' Matlab procedures
for the Prediction Error Method (see i.e.\ \cite{SODSTO} or \cite{LJUNG})
and for the N4SID method (see \cite{VODM} or \cite{LJUNG}):
The former is a multivariable extension of the classical
approach to ARMAX identification, while the latter is a standard algorithm
in the modern field of subspace identification.
In order
to obtain a comparison reasonably independent of the specific data set,
we have performed $50$ independent runs each with $100$ samples of $y$.
In such a way we have obtained $50$ different estimates $\hat{\Phi}_{{\rm M},i}$,
${\rm M=H, ME, PEM, N4SID}, \;i=1,2,\dots,50$, for each method.
%

We have then defined
\begin{equation}
E_{\rm H}(\vartheta) := \frac{1}{50} \sum_{i=1}^{50} \|\hat{\Phi}_{{\rm H},i}({\rm e}^{j\vartheta})- \Phi ({\rm e}^{j\vartheta})\|,
\end{equation}
where $\|\cdot\|$ denotes the spectral norm. This  is understood as the average estimation error of our method at each frequency.
Similarly, we have defined the average errors $E_{\rm ME}(\vartheta)$, $E_{\rm PEM}(\vartheta)$,  and $E_{\rm N4SID}(\vartheta)$ of the other methods.
In the each of the  plots  of Figure \ref{COMPARISONFIGURE}, we depict the average error of our method $E_{\rm H}(\vartheta) $ together with the average error of one of the other methods.
More explicitly, the first diagram shows the error for the Hellinger
approximation method and for the maximum entropy spectrum described in
\cite{G5}. The second diagram shows the error for the Hellinger
approximation and for the spectrum obtained via MATLAB's PEM
identification method. The third diagram shows the same for Hellinger
approximation and MATLAB's N4SID method. The Hellinger approximation based approach appears to perform better or much better than the other methods. 
The simulation yields similar results with $N=200$ data points.
With $N=300$ data samples, PEM and N4SID perform as well as our method.

Of course, one should always take into account the complexity of the resulting spectrum.
In this example,
$G(z)$ being of order $9$, the resulting spectral factor (or ``model'')
produced by the Hellinger approximation
has order $18$, whereas the corresponding maximum entropy model has order $9$ and
both N4SID and PEM usually choose order $10$.

In our simulation, the norm of the difference of two estimates
produced by PEM or by N4SID is sometimes very large when compared to
the norm of the difference between any two of the estimates produced by our method.
That is, although PEM and N4SID are provably consistent as $N\rightarrow\infty$,
when few data are available both of them may introduce occasional artifacts,
which are well visible as ``peaks'' in figure \ref{COMPARISONFIGURE} (a ``peak'' 
in the 50-run
average is due to a very high error in one of the runs, not to a systematic error).
Our method appears to be more robust in this respect.

\section{Conclusion}
In this paper, we considered the the new approach to multivariate spectrum approximation
problem with respect to the multivariable Hellinger distance, which was proposed in \cite{FPRMULTIV}.
We developed in detail the matricial Newton algorithm which was sketched there, and proved
its global convergence. Finally, we  described an application of this approach to
spectral estimation, and tested it against the well-known PEM and N4SID
algorithms.

It appears  that approximation
in the Hellinger distance may be a useful tool to gain insight into
the dynamics of a multivariate process when fewer data are available. 
In particular, simulations suggest that this method is less prone to
produce artifacts than PEM and N4SID.
Another advantage of our method and of the maximum entropy paradigm
is that
a higher resolution estimate in a prescribed frequency band
can 
be easily achieved by  properly placing some poles of $G(z)$
close to the unit circle and with phase in the prescribed band.

{\em Numerical}
robustness of the algorithm with respect to the number and the position of the poles
is an open challenge.
Also, the analysis of the achievable precision of the results
(in a statistical sense)
has 
still to be developed.  

\section*{Acknowledgments}
The detailed comments of the anonymous reviewers are gratefully acknowledged.

\end{document}